\documentclass[12pt]{article}
\usepackage{amssymb,amsmath}

\newcommand{\Zp}{{\mathbb Z}_p}
\newcommand{\Zps}{{\mathbb Z}_{p^2}}

\newcommand{\Zl}{{\mathbb Z}_l}

\newcommand{\Qp}{{\mathbb Q}_p}
\newcommand{\Fp}{{\mathbb F}_{\!p}}
\newcommand{\Fpb}{\overline{\mathbb F}_{\!p}}
\newcommand{\rhobar}{\overline{\rho}}
\newcommand{\sigmabar}{\overline{\sigma}}

\newcommand{\xbar}{\overline{x}}

\renewcommand{\hbar}{\overline{h}}
\newcommand{\uG}{\underline{G}}
\newcommand{\xhat}{\hat{x}}
\newcommand{\Fps}{{\mathbb F}_{\!p^2}}

\newcommand{\Q}{{\mathbb Q}}

\newcommand{\univ}{(\underline{\hat{A}},\underline{\hat{\imath}})}
\newcommand{\univZ}{(\underline{A},\underline{i},\underline{Z})}
\newcommand{\p}{{\mathbb P}}

\newcommand{\Z}{{\mathbb Z}}
\newcommand{\M}{{\mathbb M}}
\newcommand{\GL}{\mathbb{GL}}
\newcommand{\SL}{\mathbb{SL}}
\newcommand{\C}{{\mathbb C}}
\newcommand{\ra}{\rightarrow}
\newcommand{\lra}{\longrightarrow}
\renewcommand{\O}{{\cal O}}
\newcommand{\ODel}{{\cal O}_{\!1,N^-}}
\newcommand{\ODelN}{{\cal O}_{\!N^+,N^-}}

\newcommand{\X}{{\cal X}}
\renewcommand{\H}{{\cal H}}

\newcommand{\s}{{\cal S}}
\newcommand{\T}{{\cal T}}
\newcommand{\E}{{\cal E}}

\newcommand{\R}{{\cal R}}
\newcommand{\F}{{\cal F}}
\newcommand{\G}{{\cal G}}
\newcommand{\I}{{\cal I}}

\newcommand{\m}{{\cal M}}
\newcommand{\q}{{\cal Q}}
\renewcommand{\P}{{\cal P}}
\renewcommand{\c}{{\cal C}}
\newcommand{\hatU}{\,\hat{\cal U}}

\newcommand{\Spf}{\mbox{\rm Spf}}
\newcommand{\Spec}{\mbox{\rm Spec}}
\newcommand{\End}{\mbox{\rm End}}

\newcommand{\Lie}{\mbox{\rm Lie}}
\newcommand{\Aut}{\mbox{\rm Aut}}
\newcommand{\Tr}{\mbox{\rm Tr}}
\newcommand{\Nr}{\mbox{\rm Nr}}
\newcommand{\dst}{\displaystyle}
\newcommand{\proof}{\noindent {\em Proof: }}
\newcommand{\qed}{\hspace{\fill}$\square$\medskip}
\newcommand{\Qed}{\hspace{\fill}$\square$}

\newtheorem{theorem}{Theorem}[section]
\newtheorem{lemma}[theorem]{Lemma}
\newtheorem{prop}[theorem]{Proposition}
\newtheorem{cor}[theorem]{Corollary}
\newtheorem{definition}[theorem]{Definition}
\newenvironment{remark}{\noindent\refstepcounter{theorem}{\bf Remark \thetheorem} }{}
\title{Intersection numbers of Heegner
divisors on Shimura curves}
\author{Kevin Keating \\
Department of Mathematics \\
University of Florida \\
Gainesville, FL 32611 \\
USA \\[.2cm]
David P. Roberts \\
Division of Science and Mathematics \\
University of Minnesota \\
Morris, MN 56267 \\
USA \\[.2cm]
{\tt keating@math.ufl.edu} \\
{\tt roberts@morris.umn.edu}}
\date{}
\begin{document}
\maketitle

     In foundational papers, 
Gross, Zagier, and Kohnen established two formulas for
arithmetic intersection numbers of certain Heegner
divisors on integral models of modular curves.  
In \cite{gz1}, only one imaginary quadratic discriminant
plays a role.  In \cite{gz2} and \cite{gkz}, two
quadratic discriminants play a role.  In this
paper we generalize the two-discriminant formula
from the modular curves $X_0(N)$ to certain Shimura curves 
defined over $\Q$.  

     Our intersection formula was stated in \cite{rob},
but the proof was only outlined there.
Independently, the general formula was given, in a 
weaker and less explicit form, in  \cite{kea}; there it was
proved completely.  This paper is thus a synthesis of
parts of \cite{rob} and \cite{kea}.
The intersection multiplicities computed here were
used in \cite{kud} to derive a relation between height
pairings and special values
of the derivatives of certain Eisenstein series. 
We note also that Zhang \cite{zhang} has generalized
all of \cite{gz1} from ground field $\Q$ to 
general totally real ground fields $F$, working
with general Shimura curves.  So we certainly
expect that all of \cite{gkz} should generalize 
similarly.  Our work here can be viewed as
a step in this direction. 

   We are happy to thank B.\ Gross and S.\ Kudla for their
encouragement during the long period in which this work was done.

\section{Eichler orders} \label{eichler}
\setcounter{equation}{0}

     Let $\Delta$ be a quaternion algebra
over $\Q$ and let $p$ be a prime.  We
begin by defining Eichler orders in
$\Delta\otimes\Qp$.  Let $\E_p$ be an order
in $\Delta\otimes\Qp$ and let $p^e$ be the
reduced discriminant of $\E_p$.  We say
that $\E_p$ is an Eichler order of type
$(p^e,1)$ if $\E_p$ contains a subring
isomorphic to $\Zp\oplus\Zp$.  We say
that $\E_p$ is an Eichler order of type
$(1,p^e)$ if $\E_p$ contains a subring
isomorphic to $\Zps$, the ring of
integers in the quadratic unramified
extension of $\Qp$.  It is easily seen
that two of these local Eichler
orders are conjugate in $\Delta\otimes\Qp$
if and only if they are of the same
type.  In fact, if $\E_p$ is an
Eichler order of type $(p^e,1)$ then
$\Delta\otimes\Q_p\cong\M_2(\Q_p)$, and
$\E_p$ is conjugate to the standard Eichler order
\begin{equation} \label{standard}
\hat{\O}_{p^e,1}=\left\{\begin{bmatrix}
a&b\\c&d\end{bmatrix}\in\M_2(\Zp):p^e\mid c\right\}.
\end{equation}
Let $\m_p$ be a maximal order in $\Delta\otimes\Qp$, and let
$\psi:\Zps\ra\m_p$ be an embedding.  If $\E_p$ is an
Eichler order of type $(1,p^e)$ then $\E_p$ is
conjugate to the standard Eichler order
\begin{equation} \label{standard2}
\hat{\O}_{1,p^e}=\psi(\Zps)+p^{f}\m_p,
\end{equation}
where $f=\lfloor e/2\rfloor$.  Note that if $e$ is even then
$\Delta\otimes\Qp\cong\M_2(\Qp)$, while if $e$ is odd then
$\Delta\otimes\Qp$ is a division ring.  For future use we
let $\hat{B}_p$ denote the quaternion division ring whose
maximal order is $\hat{\O}_{1,p}$.

     We also need to define the
notion of an orientation on a local
Eichler order $\E_p$, which we assume is
not of type $(1,1)$.  For $e\ge1$ define rings
$R_{p^e,1}=(\Z/p^e\Z)\oplus(\Z/p^e\Z)$ and
$R_{1,p^e}=\Zps/p^e\Zps$.
If $\E_p$ is an Eichler order of type $(p^e,1)$
then an orientation on $\E_p$
is defined to be a ring homomorphism
$\psi_p:\E_p\ra R_{p^e,1}$.  If
$\E_p$ is an Eichler order of type $(1,p^e)$ then an
orientation on $\E_p$ is defined to be a
ring homomorphism $\psi_p:\E_p\ra R_{1,p^e}$.
Thus if $\E_p$ is an Eichler order in $\Delta\otimes\Q_p$ which is not
isomorphic to $\M_2(\Zp)$ there
are exactly two orientations on
$\E_p$.  The usefulness of giving
orientations to our Eichler
orders may be summarized in the
statement that the automorphisms of
the oriented order $(\E_p,\psi_p)$ are
precisely the maps given by
conjugation by elements of $\E_p^{\times}$.

     To define global Eichler orders we
let $N^+=\prod p^{n_p^+}$ and $N^-=\prod
p^{n_p^-}$ be relatively prime positive
integers and set $N=N^+N^-$.  We say that an order $\E$ in
$\Delta$ is an Eichler order of type $(N^+,N^-)$ if $\E\otimes\Zp$
is an Eichler order of type $(p^{n_p^+},p^{n_p^-})$ for
every prime $p$.  An
orientation on $\E$ consists of a collection
$\{\psi_p\}_{p\mid N}$ of orientations on $\E\otimes\Zp$ for
every prime $p$ which divides $N$.

\begin{prop} \label{finite}
Let $\Delta$ be a quaternion algebra over $\Q$ and let $(N^+,N^-)$
be relatively prime positive integers. \\[\smallskipamount]
(a) $\Delta$ contains an Eichler order of type $(N^+,N^-)$ if
and only if $v_p(N^-)$ is odd precisely for those primes $p$
which are ramified in $\Delta$. \\[\smallskipamount]
(b) $\Delta$ contains only
finitely many isomorphism classes of oriented
Eichler orders of type $(N^+,N^-)$. \\[\smallskipamount]
(c) If $\Delta$ is indefinite then $\Delta$ contains at most
one isomorphism class of oriented
Eichler orders of type $(N^+,N^-)$.
\end{prop}

\proof (a) It follows from the definitions that if $\Delta$
contains an Eichler order
of type $(N^+,N^-)$ then $v_p(N^-)$ is odd if and only if
$\Delta$ is ramified at $p$.  On the other hand, if $N^-$ satisfies
this condition then one can easily construct an oriented
Eichler order of type $(N^+,N^-)$ from a maximal order in $\Delta$.
\\[\smallskipamount]
(b) Let $(\E,\{\psi_p\}_{p\mid N})$ be an oriented
Eichler order in $\Delta$ of type $(N^+,N^-)$,
let $\hat{\E}=\E\otimes\hat{\Z}$ be the profinite
completion of $\E$, and let
$\hat{\Delta}=\hat{\E}\otimes\Q$
be the ring of finite ad\`eles of $\Delta$.
Associated to each $\beta=(\beta_p)\in
\hat{\Delta}^{\times}$  there is a unique
lattice $L_{\beta}$ in $\Delta$ such that
$L_{\beta}\otimes\Zp=\beta_p(\E\otimes\Zp)$
for all primes $p$.  There is a bijection between the double
coset space $S=\Delta^{\times}\backslash
\hat{\Delta}^{\times}/\hat{\E}^{\times}$ and
the set of all isomorphism classes of
oriented Eichler orders of type
$(N^+,N^-)$, which associates to
$\beta\in\hat{\Delta}^{\times}$ the pair
$(\E_{\beta},\{\phi_p^{\beta}\}_{p\mid N})$, where
$\E_{\beta}$ is the left order of $L_{\beta}$ and
$\phi_p^{\beta}(x)=\psi_p(\beta_p^{-1} x\beta_p)$.  Since $S$
is finite \cite[III, Cor.\,5.5]{vig}, the
claim follows. \\[\smallskipamount]
(c) By the strong approximation
theorem \cite[III, Th.\,4.3]{vig} the reduced norm on
$\hat{\Delta}^{\times}$ induces a bijection between $S$ and the
set $T=\Q^{\times}\backslash\hat{\Q}^{\times}/
\Nr(\hat{\E}^{\times})$.  In fact
$\Nr(\hat{\E}^{\times})=\hat{\Z}^{\times}$,
so $T$ has just one element. \Qed

\begin{cor} \label{2order}
Let $\Delta$ be an indefinite quaternion algebra over $\Q$
ramified at the primes $p_1,\ldots,p_s$.  Set
$N^-=p_1p_2\cdots p_s$, and let $N^+$ be a positive integer
which is relatively prime to $N^-$.  Then there are Eichler
orders $\ODel\supset\ODelN$ in $\Delta$ of types $(1,N^-)$
and $(N^+,N^-)$ such that $\ODel/\ODelN$
is a cyclic group of order $N^+$.  The pair $(\ODel,\ODelN)$
is uniquely determined up to conjugacy in $\Delta$.
\end{cor}

\proof The existence of the pair $(\ODel,\ODelN)$ is clear;
what must be proved is that all such pairs are conjugate in
$\Delta$.  Since $\Delta$ is indefinite, it follows from
Proposition~\ref{finite}(c) that $\ODel$ is determined
uniquely up to conjugation.  Let $\Sigma$ be the set of
Eichler orders $\E$ of type $(N^+,N^-)$ such that
$\E\subset\ODel$ and $\ODel/\E$ is cyclic of order $N^+$.
For each prime $p$ such that $p\mid N^+$ let $\Sigma_p$
denote the set of local Eichler orders $\E_p$ of type
$(p^{n^+},1)$ such that $\E_p\subset\M_2(\Z_p)$ and
$\M_2(\Z_p)/\E_p$ is cyclic.  Then
$\SL_2(\Z_p)$ acts transitively by conjugation on $\Sigma_p$.
Therefore by the strong approximation theorem the
group of elements of $\ODel^{\times}$ with reduced norm 1 acts
transitively by conjugation on $\Sigma$.  It follows that the pair
$(\ODel,\ODelN)$ is determined uniquely up to conjugation in
$\Delta$. \medskip \qed

     Let $D_1,D_2$ be negative integers which are squares
(mod 4) such that $\Q(\sqrt{D_1})\not\cong\Q(\sqrt{D_2})$.
Let $n$ be an integer such that
$n\equiv D_1D_2\pmod2$ and let $B_n$ be the Clifford algebra
of the binary quadratic form $q_n(x,y)=D_1x^2+2nxy+D_2y^2$.
Thus $B_n$ is a quaternion algebra over $\Q$ which is
generated by elements $e_1,e_2$ such that $e_j^2=D_j$ for
$j=1,2$ and $e_1e_2+e_2e_1=2n$.  Let $g_j=(D_j+e_j)/2$ and
let $S_n=\Z[g_1,g_2]$ be the subring of $B_n$ generated by
$g_1$ and $g_2$.  Then
\begin{align}
S_n&=\Z+\Z g_1+\Z g_2+\Z g_1 g_2
\end{align}
is an order in $B_n$ with reduced discriminant
$\delta_n=(n^2-D_1D_2)/4$.  We may view $S_n$ with the reduced
norm form $\Nr$ as a quadratic space over $\Z$.  By
restricting $\Nr$ to $L_n=\Z+\Z g_1+\Z g_2$ we
get a quadratic form
\begin{align} \label{Qn}
Q_n(x,y,z)&=\Nr(x+yg_1+zg_2) \\[.1cm]
&=\dst x^2+\frac{D_1^2-D_1}{4}y^2+\frac{D_2^2-D_2}{4}z^2+D_1xy+
D_2xz+\frac{D_1D_2-n}{2}yz \label{Qnexp}
\end{align}
with determinant $2\delta_n$.

     Assume now that  $n^2<D_1D_2$ and $\gcd(D_1,D_2)=1$. 
We factor the positive integer $-\delta_n$ into relatively prime
factors $\delta_n^+$, $\delta_n^-$ using the criterion
\begin{align} \label{crit1}
p\mid \delta_n^+&\mbox{ if }\left(\frac{D_j}{p}\right)=+1
\mbox{ for at least one }j=1,2, \\[.1cm]
p\mid \delta_n^-&\mbox{ if }\left(\frac{D_j}{p}\right)=-1
\mbox{ for at least one }j=1,2, \label{crit2}
\end{align}
where $\left(\frac{D_j}{p}\right)$ is the Kronecker symbol.
Suppose $p\mid \delta_n$ and $p\nmid D_1D_2$.  Then
$D_1D_2\equiv n^2\pmod{4p}$, and hence
$\left(\frac{D_1}{p}\right)=\left(\frac{D_2}{p}\right)$.
Thus (\ref{crit1}) and (\ref{crit2})
uniquely determine the factorization $-\delta_n=\delta_n^+\delta_n^-$.
For each prime $p$ we have $p\nmid D_j$ for at least one
$j\in\{1,2\}$.
Hence $\O_{D_j}\otimes\Zp$ is isomorphic to either
$\Zp\oplus\Zp$ or $\Zps$.  In particular, if $p\mid \delta_n^+$
then $\O_{D_j}\otimes\Zp\cong\Zp\oplus\Zp$, and if $p\mid
\delta_n^-$ then $\O_{D_j}\otimes\Zp\cong\Zps$.  Since $S_n$
contains a subring isomorphic to $\O_{D_j}$, it follows that
$S_n$ is an Eichler order of type $(\delta_n^+,\delta_n^-)$.

\section{Heegner divisors on Shimura curves} \label{heeg}
\setcounter{equation}{0}

     In this section we construct a scheme $\X=\X_{N^+,N^-,m}$
associated
to the triple $(N^+,N^-,m)$ for certain values of $m$.
The scheme $\X$ is an integral model for a Shimura curve
$X$ which is defined over $\Q$.  We also define
Heegner divisors $P_D$ on $X$ and $\P_D$ on $\X$,
where $D$ is the discriminant of an order in an imaginary
quadratic field.

     The scheme $\X$ will be constructed as a moduli space
for abelian surfaces $A$ with additional structure.  Part of
the additional structure on $A$ is a ``special'' $\ODel$-action,
as defined by Drinfeld \cite[\S2A]{cd}.
Let $R$ be a ring, let $A$ be an abelian surface over
$R$, and let $i:\ODel\ra\End_R(A)$ be an embedding.  For
$a\in\ODel$ let $\Tr(a)\in\Z$ denote the reduced trace
of $a$, and let $\tau(a)$ denote the image of
$\Tr(a)$ under the natural map $\Z\ra R$.  The embedding
$i:\ODel\ra\End_R(A)$ is said to be {\em special} if the
trace of the action of $i(a)$ on $\Lie(A)$ is equal to
$\tau(a)$ for all $a\in\ODel$.  (If all the primes
$p_1,\ldots,p_s$ which ramify in $\Delta$ are invertible in
$R$ then every embedding is special.)  More generally, if
$Y$ is a scheme and $A/Y$ is an abelian surface, we
say that the embedding $i:\ODel\ra\End_Y(A)$ is special if
the induced map
\begin{equation}
i_R:\ODel\lra\End_R(A\times_Y\Spec(R))
\end{equation}
is special for every affine subscheme $\Spec(R)$ of $Y$.

     We are interested in the moduli problem for isomorphism
classes of triples $(A,i,Z)$ over a scheme $Y$, where $A/Y$ is an
abelian surface, $i:\ODel\ra\End_Y(A)$ is a
special embedding, and $Z$ is a
subgroup scheme of $A$ which is cyclic of order $N^+$ in the
sense of \cite[1.4]{km}.  Since this moduli
problem is not representable, we will add a level-$m$
structure to the problem for an appropriate value of $m$.
The choice of $m$ depends on $N=N^+N^-$ and on the
imaginary quadratic discriminants $D_1,D_2$ which will be
introduced in \S\ref{formulas}.  Let $m$
be a positive integer satisfying the following conditions:
\begin{equation} \label{cond}
\begin{array}{l}
\text{C1: $\gcd(m,N)=1$.} \\[.2cm]
\text{C2: $m=m_1m_2$ for some $m_1,m_2\ge4$ such that
$\gcd(m_1,m_2)=1$.}
\\[.2cm]
\text{C3: $p>D_1D_2/4N$ for every prime $p$ which divides $m$.}
\end{array}
\end{equation}
Fix a ring isomorphism $\ODel/m\ODel\cong\M_2(\Z/m\Z)$.
Then $i:\ODel\ra\End_Y(A)$ induces a ring homomorphism
\begin{equation}
i_m:\M_2(\Z/m\Z)\lra\End_Y(A[m]),
\end{equation}
where $A[m]$ denotes the $m$-torsion subgroup scheme of $A$.
A $\Gamma_1(m)$-structure on the pair $(A,i)$ is defined to
be a point $\beta$ in the kernel of
$i_m\left(\begin{bmatrix}1&0\\0&0\end{bmatrix}\right)$ which
is defined over $Y$ and has exact order $m$ in the sense of
\cite[1.4]{km}.

     Let $\F_m$ denote the functor from schemes to sets which
associates to a scheme $Y$ the set of isomorphism classes of
4-tuples $(A,i,Z,\beta)$, where
\begin{enumerate}
\item $A$ is an abelian surface defined over $Y$.
\item $i:\ODel\ra\End_Y(A)$ is a special embedding.
\item $Z$ is a cyclic subgroup scheme of $A$ of order $N^+$
which is defined over $Y$ and stabilized by $i(\ODelN)$.
\item $\beta$ is a $\Gamma_1(m)$-structure on $(A,i)$.
\end{enumerate}
It follows from \cite[Prop.\,4.4]{cd} and \cite[Lemma 2.2]{bu}
that for $j=1,2$ the restriction of $\F_{m_j}$ to
$\Z[1/m_j]$-schemes is represented by a scheme
over $\Z[1/m_j]$.  Therefore by \cite[4.3.4]{km} the
restriction of $\F_m$ to $\Z[1/m_j]$-schemes is also
represented by a $\Z[1/m_j]$-scheme.  It follows
that the functor $\F_m$ is represented by a scheme $\X$.
By \cite[Prop.\,4.4]{cd}, the scheme $\X\otimes\Z[1/m]$ is
projective over $\Z[1/m]$.   Hence $X:=\X\otimes\Q$ is a
projective curve over $\Q$.

     Let $k$ be a field, let $x\in\X(k)$, and let
$(A_x,i_x,Z_x,\beta_x)$ be the 4-tuple which corresponds to
$x$.  An
endomorphism of the triple $(A_x,i_x,Z_x)$ is defined to be
an endomorphism of $A_x$ which stabilizes $Z_x$ and commutes
with $i_x(a)$ for every $a\in\ODel$.  Let $D$ be a negative
integer which is a square (mod 4), let $K=\Q(\sqrt{D})$, and
let $\O_D=\Z[(D+\sqrt{D})/2]$ be the order of discriminant
$D$ in $K$.  It follows from Proposition~\ref{product} below
that there are only finitely many $x\in X(\C)$ such that
$\End(A_x,i_x,Z_x)\cong\O_D$.  Therefore we may
define a divisor $Q_D$ on $X$ by setting $Q_D=\sum\:(x)$,
where the sum is taken over all such $x$.
It follows from the definition that $Q_D$ is defined over $\Q$.  
Write $D=c^2D_0$, where $c$ is the conductor of $\O_D$ and
$D_0$ is the discriminant of $K$.  Define a divisor $P_D$ on
$X$ by setting $P_D=\sum_{b\mid c}Q_{b^2D_0}$.  Then we have
$P_D=\sum\:(x)$, where the sum is taken over points $x\in X(\C)$
such that $\O_D$ embeds as a subring in $\End(A_x,i_x,Z_x)$.
We call $P_D$ the Heegner divisor of discriminant $D$ on $X$.
Since $P_D$ is defined over $\Q$, we can also express $P_D$
as a formal sum $P_D=\sum\;(y)$ of irreducible subschemes $y$
of $X$.  Let $\P_D$ be the divisor on $\X$ obtained by
replacing each subscheme in this sum by its
closure $\overline{y}$ in $\X$.

     The following proposition gives a stringent condition
that $A$ must satisfy if $(A,i,Z)$ corresponds to a point
in the support of $Q_D$.

\begin{prop} \label{product}
Let $x\in X(\C)$ be a point in the support of $Q_D$.  Let
$\R$ be the smallest order in $K=\Q(\sqrt{D})$ which contains
$\O_D$ and whose conductor is not divisible by any prime
which is ramified in $\Delta$.  Then over $\C$ we have
$A_x\cong E_1\times E_2$, where $E_1$ and $E_2$ are elliptic
curves such that $\End(E_1)\cong\End(E_2)\cong \R$.
\end{prop}

\proof Since $\End(A_x,i_x,Z_x)\cong\O_D$ there is an embedding
of $\Delta\otimes K$ into $\End(A_x)\otimes\Q$.  Therefore by
\cite[Prop.\,6.1]{end} we see that $K$ splits $\Delta$ and that
\begin{equation}
\End(A_x)\otimes\Q\cong\Delta\otimes K \cong\M_2(K).
\end{equation}
It follows that $\End(A_x)$ is isomorphic to an order $S$
in $\M_2(K)$, and that the complex points of $A_x$ may be
identified with a quotient $\C^2/L$,
where $L$ is a $\Z$-lattice in $K^2\subset
\C^2$.  The stabilizer of $L$ in $\M_2(K)$ is $S$,
and hence for each prime $p$ the stabilizer
of $L\otimes\Zp$ in $\M_2(K\otimes\Qp)$ is $S\otimes\Zp$.

     The homomorphism $i_x:\ODel\ra\End(A_x)\cong S$ induces
a map
\begin{equation}
i_x\otimes\Zp:\ODel\otimes\Zp\lra S\otimes\Zp.
\end{equation}
If $p\nmid N^-$ then $\ODel\otimes\Zp\cong\M_2(\Zp)$, and hence
$S\otimes\Zp$ is isomorphic to an
order in $\M_2(K\otimes\Qp)$ which
contains $\M_2(\Zp)$.  Such an order
must be isomorphic to $\M_2(\R_p)$
for some order $\R_p$ in $K\otimes\Qp$.
Since $\End(A_x,i_x,Z_x)\cong\O_D$ we
have $\R_p\cong\O_D\otimes\Zp$.

     If $p\mid N^-$ then $p$ is ramified
in $\Delta$, and hence $\ODel\otimes\Zp\cong\hat{\O}_{1,p}$.
Since $K$ splits $\Delta$ we see that $K_p=K\otimes\Qp$
is a field which is a quadratic extension
of $\Qp$.  Let $\O_{K_p}=\O_K\otimes\Zp$ be the
ring of integers in $K_p$.  We will show
that $S\otimes\Zp\cong\M_2(\O_{K_p})$.
Choose a $\Qp$-embedding $\psi:K_p\ra\Delta\otimes\Qp$.  We
give $\Delta\otimes\Qp$ the structure of a
$K_p$-vector space by setting
$a\cdot v=v\psi(a)$ for $a\in K_p$,
$v\in\Delta\otimes\Qp$.  Left
multiplication gives a representation of
$\Delta\otimes\Qp$ on this 2-dimensional
$K_p$-vector space.  On the other hand,
since $S\otimes\Zp$ is isomorphic to an
order in $\M_2(K_p)$, the map $i_x\otimes\Qp$
induces a representation of
$\Delta\otimes\Qp$ on $K_p^2$.  By the Skolem-Noether
theorem these two representations are isomorphic.  Let
$\Phi:K_p^2\ra\Delta\otimes\Qp$ be a
$(\Delta\otimes\Qp)$-equivariant
isomorphism of $K_p$-vector spaces.  Since
$\ODel\otimes\Z_p\cong\hat{\O}_{1,p}$
stabilizes $L\otimes\Zp$, it stabilizes
$\Phi(L\otimes\Zp)$ as well.  Therefore
$\Phi(L\otimes\Zp)$ is a left
$(\ODel\otimes\Z_p)$-ideal, and hence also a right
$(\ODel\otimes\Z_p)$-ideal.  Since
$\psi(\O_{K_p})\subset\ODel\otimes\Z_p$ this implies
that $\Phi(L\otimes\Zp)$ is an $\O_{K_p}$-module.  Since
$\Phi(L\otimes\Zp)$ is free of rank 4 over $\Z_p$, it is free
of rank 2 over $\O_{K_p}$.  Therefore $L\otimes\Zp$ is
also a free $\O_{K_p}$-module of rank 2.
We conclude that the stabilizer $S\otimes\Zp$ of
$L\otimes\Zp$ is isomorphic to $\M_2(\O_{K_p})$.

     So far we have proved that $\End(A_x)$
is isomorphic to an order $S$ in $\M_2(K)$
such that $S\otimes\Zp\cong\M_2(\R_p)$
for all $p$, where $\R_p=\O_{K_p}$ if $p$ is ramified in
$\Delta$ and $\R_p=\O_D\otimes\Z_p$ if $p$ is not ramified in
$\Delta$.  Hence the order $\R$ in the statement of the
theorem is the unique order in $K$ such that $\R\otimes\Zp=\R_p$
for all $p$.  To complete the proof we will
show that $S$ contains a nontrivial idempotent.

     Let $L'\supset L$ be the $\O_K$-lattice generated by
$L$.  By choosing a new $K$-basis for
$K^2\subset\C^2$ we may assume that
$L'=\O_K\oplus\,\I$ for some ideal
$\I\subset\O_K$.  The ideal $\I$ may
be chosen to be relatively prime to
every $p$ such that $\R_p$ is not the
maximal order in $K_p$.  There
is an abelian surface $A'$ over $\C$
such that $A'(\C)\cong\C^2/L'\cong(\C/\O_K)\times(\C/\I)$.
The endomorphism ring of $A'$ is
\begin{align}
S'&\cong\End_{\O_K}(L') \\[.1cm]
&\cong\left\{\left[\begin{array}{cc}
a&b\\c&d\end{array}\right]:a,d\in\O_K,
\,b\in\I^{-1},\,c\in\I\right\},
\end{align}
which is a maximal order in $\M_2(K)$
containing $S$.  There is an action of
$\ODel$ on $A'$ given by the map
$i':\ODel\ra S'$ which is the
composition of the inclusion
$S\hookrightarrow S'$ with $i_x:\ODel\ra S$.  

     The inclusion $L\hookrightarrow L'$
induces an $\ODel$-equivariant isogeny $\pi:A_x\ra A'$.
The kernel of $\pi$ is a finite
subgroup $G\cong L'/L$ of $A_x(\C)$
which is stabilized by $i(\ODel)$.  Let
$G\cong\oplus_p G_p$ be the
decomposition of $G$ into its
$p$-primary components.  Then $G_p\cong
(L'\otimes\Zp)/(L\otimes\Zp)$, and
$G_p=\{0\}$ for all $p$ such that
$\R_p=\O_K\otimes\Zp$.  In particular,
$G_p=\{0\}$ if $p$ is ramified in
$\Delta$.  

     Let $p$ be a prime such that $G_p\not=\{0\}$.  Since
$S\otimes\Zp\cong\M_2(\R_p)$ is the
stabilizer of $L\otimes\Zp$ in $\M_2(K\otimes\Qp)$,
we see that $L\otimes\Zp$ is free of rank 2 over $\R_p$.  Let
$C_p\in\M_2(K\otimes\Qp)$ be a matrix
whose columns are $\R_p$-generators for $L\otimes\Zp$.
The columns of $C_p$ also serve
as $(\O_K\otimes\Zp)$-generators for $L'\otimes\Zp$,
and by the assumption on $\I$ we have
$L'\otimes\Zp=(\O_K\otimes\Zp)^2$.
Therefore $C_p\in\GL_2(\O_K\otimes\Zp)$.
Let $J=\left[\begin{array}{cc}1&0\\0&0\end{array}\right]$.
By multiplying one of the columns of $C_p$ by $1/\det(C_p)$
we get a matrix $C_p'\in\SL_2(\O_K\otimes\Zp)$ such that
$C_p'JC_p'^{-1}=C_pJC_p^{-1}$ is a nontrivial idempotent
which lies in $\End_{\R_p}(L\otimes\Zp)=S\otimes\Zp$.

     For each $p$ such that $G_p\not=\{0\}$ choose $n_p\geq1$
such that $p^{n_p}$ kills $G_p$.  By the strong approximation
theorem there exists a matrix $C\in\SL_2(\O_K)$ such that
\begin{alignat}{2}
C&\equiv C_p'&&\pmod{p^{n_p}}\mbox{ for all $p$ such that }
G_p\not=\{0\}, \\[.1cm]
C&\equiv I_2&&\pmod{\I}.
\end{alignat}
Then $e=CJC^{-1}$ is a nontrivial idempotent in $S$.  Set
$E_1=eA_x$ and $E_2=(1-e)A_x$.  Then $A_x\cong E_1\times E_2$
with $\End(E_1)\cong\End(E_2)\cong \R$.~\qed

     Let $c$ be the conductor of $\O_D$. \medskip

\begin{remark} \label{ray}
Recall that the ray
class field $K_c$ of $K$ with conductor $c\O_K$ is the
maximum abelian extension of $K$ whose ramification
conductor divides $c\O_K$.
The elliptic curves $E_1$ and $E_2$ are defined over
$K_c$, and hence the triple $(A_x,i_x,Z_x)$ is defined over
$K_{cN^+}$.  Therefore $x$ is rational over $K_{cN^+}$.
\medskip
\end{remark}

\begin{remark} \label{conductor}
If $\gcd(c,N^-)\not=1$ then by Proposition~\ref{product} we
get $Q_D=0$.  Therefore if $p$ divides $\gcd(c,N^-)$ then
$P_D=P_{D/p^2}$.  Hence we may assume
without loss of generality that $c$ is relatively prime to $N^-$.
\medskip
\end{remark}

\begin{remark} \label{inert}
Suppose $p^t\mid N^+$.  Then the order in $K$ which stabilizes
$Z\subset E_1\times E_2$ has conductor divisible by $p^t$.
Hence if $p^t\nmid c$ and $p$ is inert in $K$ then $P_D=0$.
\medskip
\end{remark}

\begin{remark} \label{split}
Suppose $p\mid N^-$ and $p$ splits in $K=\Q(\sqrt{D})$.  Then
$\ODel\otimes\Z_p\cong\O_{1,p}$ cannot be embedded in
$\M_2(K)\otimes\Q_p\cong\M_2(\Q_p\oplus\Q_p)$,
so we have $P_D=0$ in this case as well.
\medskip
\end{remark}

     Associated to each point in the support of $P_D$ is a
collection of homomorphisms $\{\omega_p\}_{p\mid N}$ which
is analogous to an orientation.  The following well-known
fact will be used to construct these homomorphisms.

\begin{lemma} \label{known}
Let $R$ be a (possibly noncommutative) ring with 1 and let
$M$ be a free left $R$-module of rank 1 generated by $e\in M$.
For $\phi\in\End_R(M)$ define $f(\phi)\in R$
by the formula $\phi(e)=f(\phi)e$.
Then the map $f:\End_R(M)\ra R^{op}$ is an
isomorphism of rings, uniquely determined
by $M$ up to conjugation by units in $R^{op}$.
\end{lemma}

     Fix an orientation $\{\phi_p\}_{p\mid N}$ on $\ODelN$
and let $x\in X(\C)$ be a point in the support of $P_D$.
Let $p$ be a prime which divides $N$, let $T_p(A_x)$ be the
$p$-adic Tate module of $A_x$, and let $U_p^x\supset T_p(A_x)$
be the lattice
which corresponds to the $p$-primary part of $Z_x$.  Since
$T_p(A_x)$ is a left module over the maximal order
$\ODel\otimes\Zp$, it follows from \cite[Th.\,18.7]{mo}
that $T_p(A_x)$ is free of rank 1 over $\ODel\otimes\Zp$.
The $(\ODel\otimes\Zp)$-module structure $i_x\otimes\Zp$ on
$T_p(A_x)$ induces an $(\ODelN\otimes\Zp)$-module structure
$i_p$ on $U_p^x$.  If $p\mid N^-$ it follows that
$U_p^x=T_p(A_x)$ is free of rank 1 as a left module over
$\ODelN\otimes\Zp=\ODel\otimes\Z_p$.  If $p\mid N^+$ we may
identify $\ODel\otimes\Z_p$ with $\M_2(\Z_p)$ and
$\ODelN\otimes\Z_p$ with the standard local Eichler order
$\hat{\O}_{p^{n_p^+},1}$ defined in
(\ref{standard}).  Hence there is a generator $e$ for
$T_p(A_x)$ over $\ODel\otimes\Z_p\cong\M_2(\Z_p)$ such that
$U_p^x/T_p(A_x)$ is generated by
$\begin{bmatrix}p^{-n_p^+}&0\\0&0\end{bmatrix}\cdot e$.
It follows that $U_p^x$ is free of rank 1 over
$\ODelN\otimes\Z_p$.

     Let $p$ be a prime which divides $N$.  Since $U_p^x$ is
free of rank 1 over $\ODelN\otimes\Z_p$, by Lemma~\ref{known}
we get a ring isomorphism
\begin{equation}
\End(U_p^x,i_p)\lra\ODelN^{op}\otimes\Zp.
\end{equation}
Since
\begin{equation}
\End(U_p^x,i_p)\cong\End(A_x,i_x,Z_x)\otimes\Zp,
\end{equation}
the orientation $\phi_p$ on
$\ODelN\otimes\Zp$ induces a homomorphism
\begin{equation}
\omega_p^x:\End(A_x,i_x,Z_x)\otimes\Zp\lra R_{p^{n_p^+},p^{n_p^-}},
\end{equation}
where $n_p^+=v_p(N^+)$ and $n_p^-=v_p(N^-)$.

     For each $x$ in the support of $P_D$ there are two
embeddings of $\O_D$ into $\End(A_x,i_x,Z_x)$.  Choose one
of these and call it $\rho_x$.  For $p\mid N$, $p\not=2$
set $R_{(p)}=R_{p^{n_p^+},p^{n_p^-}}$.
By composing $\rho_x\otimes\Zp$ with $\omega_p^x$ we get a
homomorphism
$\lambda_p^x:\O_D\otimes\Zp\ra R_{(p)}$.  If $2\mid N^+$ set
$R_{(2)}=R_{2^{n_2^++1},1}$, while if $2\mid N^-$ set
$R_{(2)}=R_{1,2^{n^-+1}}=R_{1,4}$.  In either case there
is a ring homomorphism
\begin{alignat}{2}
\lambda_2^x&:1\otimes\Z_2+2\O_D\otimes\Z_2\lra R_{(2)}
\end{alignat}
defined by
\begin{equation}
\lambda_2^x(1\otimes\alpha+2\beta)=
\alpha+2\cdot\omega_2^x\circ(\rho_x\otimes\Z_2)(\beta)
\pmod{2^{v_2(N)+1}}
\end{equation}
for $\alpha\in\Z_2$, $\beta\in\O_D\otimes\Z_2$.
For every prime $p$ such that $p\mid N$ set $a_p=v_p(2N)$.
Then $\lambda_p^x(\sqrt{D})^2\equiv D\pmod{p^{a_p}}$, and the
ring homomorphism $\lambda_p^x$ is determined by the value of
$\lambda_p^x(\sqrt{D})$.  Let $\iota_p$ denote the natural
involution of the ring $R_{(p)}$.  It follows from the
definition of $\lambda_p^x$ that
$\iota_p(\lambda_p^x(\sqrt{D}))=-\lambda_p^x(\sqrt{D})$.

     Assume that the conductor $c$ of $\O_D$ is
relatively prime to $N$.  For each $p\mid N$ let $b_p$ be
an element of $R_{(p)}$ such that $b_p^2\equiv D\pmod{p^{a_p}}$
and $\iota_p(b_p)=-b_p$.  Since $p\nmid c$, there are two
possibilities for $b_p$.  Write $b=(b_p)_{p\mid N}$ and define
\begin{align}
V_b&=\{x\in\text{Supp}(P_D):\lambda_p^x(\sqrt{D})=b_p
\text{ for every }p\mid N\}.
\end{align}
Then we get a divisor $P_{D,b}=\sum_{x\in V_b}(x)$
on $X$ such that $P_D=\sum_bP_{D,b}$.
In general $P_{D,b}$ is not
defined over $\Q$ and depends on the choices of the
$\rho_x$, but the sum $P_{D,\pm b}=P_{D,b}+P_{D,-b}$ is a
well-defined divisor over $\Q$.
Define $\P_{D,\pm b}$ to be the closure of $P_{D,\pm b}$ in
$\X$.  Then $\P_{D,\pm b}$ is defined over $\Z$ and doesn't
depend on the $\rho_x$.  The divisors $P_D$ and $\P_D$
have natural sum decompositions $P_D=\sum P_{D,\pm b}$ and
$\P_D=\sum\P_{D,\pm b}$.

\section{Intersection formulas}
\label{formulas}
\setcounter{equation}{0}

     In this section we define the arithmetic intersection
number $\langle\q_1\cdot\q_2\rangle_{\X}$ of two divisors
$\q_1,\q_2$ on $\X$.  We then give formulas for computing
$\langle\P_{D_1,\pm b_1}\cdot\P_{D_2,\pm b_2}\rangle_{\X}$
and ${\langle\P_{D_1}\cdot\P_{D_2}\rangle_{\X}}$
in certain cases.

     We wish to define a $\Q$-linear pairing
$\langle\;\cdot\;\rangle_{\X}$ of divisors on $\X$ which
intersect properly on regular points of $\X$.  It suffices
to define $\langle T_1\cdot T_2\rangle_{\X}$ for dimension-1
subschemes $T_1$, $T_2$ of $\X$ whose intersection is supported
on a finite set of closed points of $\X$, each of which is
regular.  In this case we have $T_1\cap T_2\cong\Spec\,R$
for some finite ring $R$, where $T_1\cap T_2$ is understood
to mean $T_1\times_{\X}T_2$.  The arithmetic intersection number
of $T_1$ with $T_2$ is defined to be
\begin{equation} \label{component}
\langle T_1\cdot T_2\rangle_{\X}=\log\#R.
\end{equation}
This formula extends by $\Z$-linearity to give pairings of
divisors on $\X$.  

     The following proposition implies that the intersection
of $\P_{D_1}$ and $\P_{D_2}$ is supported on a finite set of
closed points of $\X$.

\begin{prop} \label{ss}
Let $t$ be a point on $\X$ which lies in the support of the
intersection of $\P_{D_1}$ with $\P_{D_2}$ and let
$(A_t,i_t,Z_t,\beta_t)$ be the corresponding 4-tuple.  Then
$t$ is a closed point of characteristic $p>0$, and
over $\Fpb$ we have $A_t\cong E\times E$ for any
supersingular elliptic curve $E$.
\end{prop}

\proof Since $\Q(\sqrt{D_1})\not\cong\Q(\sqrt{D_2})$ the
images of $\O_{D_1}$ and $\O_{D_2}$ generate a subalgebra
of $\End(A_t,i_t,Z_t)$ with $\Z$-rank $\geq4$. 
Therefore by \cite[Prop.\,6.1]{end}, $t$ is a point of
characteristic $p>0$, and $A_t$ is isogenous to the product
of two supersingular elliptic curves.  Hence by
Proposition~\ref{product}, $A_t$ is actually isomorphic
to the product of two supersingular elliptic curves over $\Fpb$.
Since $t\in\xbar$ for some $x$ in the support of $P_{D_1}$,
it follows that $t$ is closed in $\X$.
A theorem of Deligne \cite[Th.\,3.5]{shioda} says that the
isomorphism class of the product of two
supersingular elliptic curves over
$\Fpb$ does not depend on the factors.
Therefore $A_t\cong E\times E$ for any
supersingular elliptic curve $E$ over $\Fpb$. \qed

     Let $D_1,D_2$ be negative integers
which are squares (mod 4) such that
$\Q(\sqrt{D_1})\not\cong\Q(\sqrt{D_2})$.
For $i=1,2$ let $c_i$ be the conductor of the order $\O_{D_i}$.
For each prime $l$ such that $l\mid N$ we make the following
additional assumptions:
\begin{align}
&\mbox{$l$ divides neither $c_1$ nor
$c_2$, and} \label{assume1} \\[.2cm]
&\mbox{$l$ divides at most one of
$D_1,D_2$.} \label{assume2}
\end{align}
Assumption (\ref{assume1}) guarantees
that $P_{D_1}$ and $P_{D_2}$ are Heegner
divisors in the sense of \cite{birch}, and in any case
may be made without loss of generality
for $l\mid N^-$ by Remark~\ref{conductor}.  Assumption
(\ref{assume2}) implies that every point in the support of
the intersection of $\P_{D_1}$ with $\P_{D_2}$ is regular
(see Corollary~\ref{regular}).
It follows that $\langle\P_{D_1}\cdot\P_{D_2}\rangle_{\X}$ is
well-defined.

     Let $p$ be a prime.  If $p$ is unramified in $\Delta$ let
$\Delta(p)$ denote the quaternion algebra over $\Q$ which is
ramified at $\infty,p,p_1,\ldots,p_s$.  If $p=p_j$ is ramified
in $\Delta$ let $\Delta(p)$ denote the quaternion algebra over
$\Q$ which
is ramified at $\infty,p_1,\ldots,p_{j-1},p_{j+1},\ldots,p_s$.
Let $\s_p$ denote the set of Eichler
orders $\E\subset\Delta(p)$ of type $(N^+,pN^-)$; it follows
from Proposition~\ref{finite}(a) that $\s_p$ is not empty.
We view the elements of $\s_p$ as lattices in
$\Delta(p)$ with the $\Z$-valued quadratic
forms induced by the reduced norm form
on $\Delta(p)$.  The assumption that $\E$ has type
$(N^+,pN^-)$ determines the isometry class
of the reduced norm form on $\E\otimes\Zl$
for every finite prime $l$, and $\E\otimes{\mathbb R}
\cong\Delta(p)\otimes{\mathbb R}$ is positive definite
since $\Delta(p)$ is ramified at $\infty$.  Therefore every
$\E\in\s_p$ belongs to the same genus of quadratic spaces.
Let $\G_p$ be a set of representatives of the proper
equivalence classes of this genus.

     Let $n$ be an integer such that $n^2<D_1D_2$ and
$n\equiv D_1D_2\pmod2$.
Let $L$ be a quadratic space over $\Z$ with finite rank and
let $w_L$ denote the number of proper self-isometries of $L$.
Define $\R_L(Q_n)$ to be the number of representations on
$L$ of the quadratic form $Q_n$ defined in (\ref{Qn}).
Also let $r$ be the number of
distinct prime divisors of $N$, and set
$\eta(m)=\frac12m^2\cdot\prod_{p\mid m}(1-p^{-2})$.

     We now state the first version of our intersection formula.

\begin{theorem} \label{repnum}
Let $D_1,D_2$ satisfy assumptions (\ref{assume1}) and
(\ref{assume2}), and assume that $m$ satisfies conditions
C1, C2, and C3 in (\ref{cond}).  Then the arithmetic
intersection number
of $\P_{D_1}$ with $\P_{D_2}$ on $\X$ is given by
\begin{equation} \label{repform}
\langle\P_{D_1}\cdot\P_{D_2}\rangle_{\X}=2^{r-1}\cdot
\eta(m)\cdot\sum_{p<\infty}
\left(\underset{n^2\equiv D_1D_2\;(4pN)}{\sum_{n^2<D_1D_2}}\left(\sum_{L\in\G_p}\;\frac{\R_L(Q_n)}{w_L}
\right)\cdot\alpha_p(Q_n)\right)\cdot\log p,
\end{equation}
where the local intersection multiplicities $\alpha_p(Q_n)$
are computed in (\ref{unram}) and (\ref{ram}).  
\end{theorem}

\begin{remark}
The inner sum on the right side of (\ref{repform})
is a representation number in the sense of Siegel (see for
instance \cite[p.\,377]{ca}).
\medskip
\end{remark}

\begin{remark} \label{eta}
Let $k$ be an algebraically closed field whose characteristic
does not divide $m$, and let $(A,i)$ be an abelian surface
with special $\ODel$-embedding defined over $k$.  Then the
pair $(A,i)$ admits $2\eta(m)$ different
$\Gamma_1(m)$-structures.
\medskip
\end{remark}

     By strengthening assumption (\ref{assume2}) we get a
formula for
$\langle\P_{D_1,\pm b_1}\cdot\P_{D_2,\pm b_2}\rangle_{\X}$
which is stated explicitly in terms of finite Dirichlet series.

\begin{definition} \label{Dirichlet}
For $p$ prime and $e\ge0$ define
\begin{align}
L_{p^e,1}(s)&=1+p^{-s}+p^{-2s}+\dots+p^{-es} \\
L_{1,p^e}(s)&=1-p^{-s}+p^{-2s}-\dots+(-1)^ep^{-es}.
\label{L1pe}
\end{align}
For relatively prime positive integers $M^+,M^-$ define
\begin{equation} \label{LMM}
L_{M^+,M^-}(s)=\prod_{p^e\parallel M^+}L_{p^e,1}(s)\cdot
\prod_{p^e\parallel M^-}L_{1,p^e}(s).
\end{equation}
\end{definition}

     For each prime $p$ such that $p\mid N$ set $a_p=v_p(2N)$.
Since $\iota_p((b_j)_p)=-(b_j)_p$ we have
$(b_1)_p(b_2)_p\equiv h_p\pmod{p^{a_p}}$ for some $h_p\in\Z$.
Let $h$ be an integer such that
\begin{alignat}{3} \label{h1}
h&\equiv h_p&&\pmod{p^{a_p}}&\qquad&\text{for all }p\mid N, \\
h&\equiv D_1D_2&&\pmod{2}&&\text{if $2\nmid N$}. \label{h2}
\end{alignat}
These congruences determine the class of $h\pmod{2N}$.

\begin{theorem} \label{explicit}
Assume that $D_1,D_2$ satisfy assumption (\ref{assume1}),
$\gcd(D_1,D_2)=1$, and $\P_{D_j,\pm b_j}\not=0$ for
$j=1,2$.  Assume further that $m$ satisfies
conditions C1, C2, and C3 in (\ref{cond}).  Then the
arithmetic intersection
number of $\P_{D_1,\pm b_1}$ with $\P_{D_2,\pm b_2}$ on
$\X$ is given by
\begin{equation} \label{explicitform}
\langle\P_{D_1,\pm b_1}\cdot\P_{D_2,\pm b_2}\rangle_{\X}=
\eta(m)\cdot\underset{n\equiv\pm h\;(2N)}{\sum_{n^2<D_1D_2}}\;
L_{\delta_n^+/N^+,\delta_n^-/N^-}'(0),
\end{equation}
where $\pm h=\pm h(\pm b_1,\pm b_2)$ is determined by
(\ref{h1}) and (\ref{h2}), and $\delta_n^+$, $\delta_n^-$ are
determined by (\ref{crit1}) and (\ref{crit2}).
\end{theorem}

\begin{remark} \label{Ndel}
If $p\mid N^+$ then since $\P_{D_j,\pm b_j}\not=0$,
Remark~\ref{split} implies that $p$ is not inert in
$\Q(\sqrt{D_j})$.  Similarly, Remark~\ref{split} implies
that if $p\mid N^-$ then $p$ is not split in
$\Q(\sqrt{D_j})$.  Since $N\mid\delta_n$ for every $n$ such
that $n\equiv\pm h\pmod{2N}$, we get $N^+\mid\delta_n^+$ and
$N^-\mid\delta_n^-$.  Hence the right side of
(\ref{explicitform}) is well-defined.
\medskip
\end{remark}

     For each value of $\pm h\pmod{2N}$ such that $h^2\equiv
D_1D_2\pmod{4N}$
there are $2^{r-1}$ pairs $(\pm b_1,\pm b_2)$ such that
$\pm b_1b_2\equiv\pm h\pmod{2N}$.  Therefore by summing
(\ref{explicitform}) over all $(\pm b_1,\pm b_2)$ 
we get a formula for
$\langle\P_{D_1}\cdot\P_{D_2}\rangle_{\X}$.

\begin{cor}
Assume that $D_1,D_2$ satisfy assumption (\ref{assume1}),
$\gcd(D_1,D_2)=1$, and $m$ satisfies conditions C1, C2, and
C3 in (\ref{cond}).  Then 
\begin{equation}
\langle\P_{D_1}\cdot\P_{D_2}\rangle_{\X}=2^{r-1}\cdot
\eta(m)\cdot\underset{n^2\equiv D_1D_2
\;(2N)}{\sum_{n^2<D_1D_2}}\;
L_{\delta_n^+/N^+,\delta_n^-/N^-}'(0).
\end{equation}
\end{cor}

\section{Intersection points}
\label{points}
\setcounter{equation}{0}

     Let $t$ be a point in the support of the intersection of
$\P_{D_1}$ with $\P_{D_2}$, and let $(A_t,i_t,Z_t,\beta_t)$ be
the corresponding 4-tuple.  In this section we study the
endomorphism ring of the triple $(A_t,i_t,Z_t)$.  In
particular, we show that $\End(A_t,i_t,Z_t)$ is an Eichler order,
and we construct an orientation on $\End(A_t,i_t,Z_t)$ which is
induced by the orientation $\{\phi_l\}_{l\mid N}$ on $\ODelN$.

     Let $T_1$, $T_2$ be dimension-1 subschemes of $\X$ whose
intersection is supported on a finite set of regular closed
points of $\X$.  Recall that the arithmetic intersection number
of $T_1$ with $T_2$ is defined to be
$\langle T_1\cdot T_2\rangle_{\X}=\log\#R$,
where $T_1\cap T_2\cong\Spec\,R$.  In practice, we will compute
$\langle T_1\cdot T_2\rangle_{\X}$ as a sum
\begin{equation} \label{pdecomp}
\langle T_1\cdot T_2\rangle_{\X}=
\sum_{p<\infty}\;(T_1\cdot T_2)_p\cdot\log p,
\end{equation}
where $(T_1\cdot T_2)_p$ is the intersection multiplicity
of $T_1$ with $T_2$ at points of characteristic $p$.  Thus
$(T_1\cdot T_2)_p$
is equal to the length of the $\Z_p$-module $R\otimes\Z_p$.
Let $W_p=W(\Fpb)$ denote the ring of Witt vectors with
coefficients in $\Fpb$.  Then $(T_1\cdot T_2)_p$ may also be
computed as the length of the $W_p$-module $R\otimes W_p$, or
as the intersection multiplicity of $T_1\otimes W_p$ with
$T_2\otimes W_p$ on $\X\otimes W_p$.

     Let $t$ be a point of characteristic $p$ in the support
of the intersection of $\P_{D_1}\otimes W_p$ with
$\P_{D_2}\otimes W_p$ on
$\X\otimes W_p$.  Then $t$ is rational over the residue
field $\Fpb$ of $W_p$ and thus may be viewed as an element of
$\X(\Fpb)$.  Let $E$ be a supersingular elliptic curve over
$\Fpb$ and set $\Lambda=\End(E)$.  Then
by Proposition~\ref{ss} we have $A_t\cong E\times E$, and
hence $\End(A_t)\cong\M_2(\Lambda)$.  We will
assume that $E$ is defined over $\Fp$; this implies that the
elements of $\End(E)$ are defined over $\Fps$.
It is well-known that $H_p=\Lambda\otimes\Q$ is a quaternion
algebra over $\Q$ which is ramified
at $p$ and $\infty$, and that
$\Lambda$ is a maximal order in $H_p$.
It follows from the definition of $\Delta(p)$ that there is
an embedding $h:\Delta(p)\ra\M_2(H_p)$ such that
$h(\Delta(p))$ is the commutant of
$i(\Delta)$ in $\M_2(H_p)$.  The endomorphism ring of
the triple $(A_t,i_t,Z_t)$ consists of those elements
of $\M_2(\Lambda)$ which commute with every element of
$i_t(\ODel)$ and stabilize $Z_t$.  Therefore $\End(A_t,i_t,Z_t)$
is identified via $h^{-1}$ with an order $\E$ in $\Delta(p)$.

     We now determine the completions of
$\End(A_t,i_t,Z_t)\cong \E$ at the finite places of $\Q$.
Let $l\not=p$, let $T_l(A_t)$ denote the
$l$-adic Tate module of $A_t$, and let $U_l^t\supset T_l(A_t)$ be
the lattice which corresponds to the $l$-primary part of $Z_t$.
As in \S\ref{heeg}, $U_l^t$ is a free
$(\ODelN\otimes\Zl)$-module of rank 1.
Using Lemma~\ref{known} we get an isomorphism
\begin{equation}
\End(A_t,i_t,Z_t)\otimes\Zl\cong\ODelN^{op}\otimes\Zl.
\end{equation}
Therefore if $l\not=p$ then $\End(A_t,i_t,Z_t)\otimes\Zl$ is a
local Eichler order of type $(l^{n_l^+},l^{n_l^-})$.

     We now consider the completion of $\End(A_t,i_t,Z_t)$
at $p$.  Leaving the subgroup $Z_t$ aside we see that
$\End(A_t,i_t)\otimes\Zp$ is the
commutant of $i(\ODel)\otimes\Zp$ in
$\End(A_t)\otimes\Zp\cong\M_2(\hat{\O}_{1,p})$.
If $p\nmid N^-$ then $i_t(\ODel)\otimes\Zp\cong\M_2(\Zp)$
stabilizes $\hat{\O}_{1,p}\oplus\hat{\O}_{1,p}$, and hence
$i_t(\ODel)\otimes\Zp=g\M_2(\Zp)g^{-1}$
for some $g\in\GL_2(\hat{\O}_{1,p})$.
Since $g\M_2(\hat{\O}_{1,p})g^{-1}=\M_2(\hat{\O}_{1,p})$ and
the commutant of $\M_2(\Z_p)$ in $\M_2(\hat{\O}_{1,p})$ is
$\hat{\O}_{1,p}\cdot I_2$, it follows that
\begin{equation} \label{O1p}
\End(A_t,i_t)\otimes\Zp\cong g(\hat{\O}_{1,p}\cdot I_2)g^{-1}
\cong\hat{\O}_{1,p}.
\end{equation}
Thus if $p\nmid N$ then
$\End(A_t,i_t,Z_t)\otimes\Zp\cong\hat{\O}_{1,p}$.
If $p\mid N^+$ we may assume that $p\nmid D_1$ by
(\ref{assume2}).  Then Remark~\ref{inert} implies that
$p$ is split in $\Q(\sqrt{D_1})$, which contradicts the
fact that $\O_{D_1}\otimes\Zp$ embeds in
$\End(A_t,i_t,Z_t)\otimes\Zp\cong\hat{\O}_{1,p}$.
So in fact assumption (\ref{assume2})
implies that there are no intersection
points in characteristic $p$ if $p\mid N^+$.

     Finally, if $p\mid N^-$ we need to consider
embeddings of $\ODel\otimes\Zp\cong\hat{\O}_{1,p}$ into
$\M_2(\hat{\O}_{1,p})$ which are induced by special
embeddings $i_t:\ODel\ra\End(A_t)$.  We can write
$\hat{\O}_{1,p}=\Zps+\Zps\pi$,
where $\Zps\subset\hat{\O}_{1,p}$ is
the ring of integers in an unramified
quadratic extension of $\Qp$, and $\pi$ is an element
of $\hat{\O}_{1,p}$ which  normalizes $\Zps$ and
satisfies $\pi^2=p$.  Consider first the
embedding of $\Zps$ into $\M_2(\hat{\O}_{1,p})$.
It is not hard to show that
since $i_t$ is special this embedding
is conjugate by an element of
$\GL_2(\hat{\O}_{1,p})$ to the map
\begin{equation} \label{zps}
x\longmapsto\left[\begin{array}{cc}x&0\\0&x'
\end{array}\right],
\end{equation}
where $x'=\pi x\pi^{-1}$ is the Galois
conjugate of $x$ over $\Qp$.
The commutant in $\M_2(\hat{\O}_{1,p})$ of
the image of (\ref{zps}) is
\begin{equation} \label{cent}
C=\left\{\left[\begin{array}{cc}a&b\pi\\
c\pi&d\end{array}\right]:a,b,c,d\in\Zps
\right\}.
\end{equation}
The image of $\pi$ under $i_t\otimes\Z_p$ is a matrix of the form
\begin{equation}
\Pi=\begin{bmatrix}e\pi&f\\g&h\pi\end{bmatrix},
\end{equation}
with $e,f,g,h\in\Zps$ and $\Pi^2=pI_2$.  If $p\nmid f$ then
$M=\begin{bmatrix}1&0\\e\pi&f\end{bmatrix}$ is a unit in
$C$, and $M\Pi M^{-1}=\begin{bmatrix}0&1\\p&0\end{bmatrix}$.
Similarly, if $p\nmid g$ there is $M\in C^{\times}$
such that $M\Pi M^{-1}=\begin{bmatrix}0&p\\1&0\end{bmatrix}$.
If $p\mid f$ and $p\mid g$ then by an iterative procedure one
constructs $M\in C^{\times}$ such that
$M\Pi M^{-1}=\begin{bmatrix}\pi&0\\0&\pi\end{bmatrix}$.
Hence, up to conjugation by units in $C$,
there are three possibilities for $(i_t\otimes\Z_p)(\pi)$,
namely
\begin{equation} \label{cases}
\Pi_1=\left[\begin{array}{cc}0&1\\p&0\end{array}\right]\hspace{1cm}
\Pi_2=\left[\begin{array}{cc}0&p\\1&0\end{array}\right]\hspace{1cm}
\Pi_3=\left[\begin{array}{cc}\pi&0\\0&\pi\end{array}\right].
\end{equation}

     Corresponding to the matrices in (\ref{cases}) are
embeddings $i_p^1$, $i_p^2$, and $i_p^3$ of $\hat{\O}_{1,p}$
into $\M_2(\hat{\O}_{1,p})$ such that $i_p^j(\pi)=\Pi_j$.  If
$i:\ODel\ra\End(E\times E)$ is a
special embedding then $i\otimes\Z_p$ is conjugate to exactly
one of these embeddings.  We say that $i$ is of type $j$ if
$i\otimes\Z_p$ is conjugate to $i_p^j$.
The commutant of the image of $i_p^1$
consists of matrices of the form (\ref{cent}) with $d=a$ and
$c=pb$; the commutant of the image of $i_p^2$
consists of matrices of the form
(\ref{cent}) with $d=a$ and $b=pc$; and
the commutant of the image of $i_p^3$
consists of matrices of the form
(\ref{cent}) with $a,b,c,d\in\Zp$.  It follows that for
$j=1,2$ the commutant of the image of $i_p^j$
in $\M_2(\hat{\O}_{1,p})$ is a
local Eichler order of type $(1,p^2)$, while the
commutant of the image of $i_p^3$ in
$\M_2(\hat{\O}_{1,p})$ is a local Eichler order of type
$(p,1)$.  By (\ref{assume2}) and Remark~\ref{split} we may
assume $\O_{D_1}\otimes\Zp\cong\Zps$.  Since
$\Zps$ cannot be embedded into a local
Eichler order of type $(p,1)$,
this implies that the third case does not occur.
Therefore if $p\mid N^-$ then $\End(A_t,i_t,Z_t)\otimes\Zp$
is a local Eichler order of type $(1,p^2)$.

     The following definition characterizes the set of
potential intersection points in characteristic $p$.

\begin{definition}
Let $p$ be a prime.
\begin{enumerate}
\item If $p\nmid N$ define $\T_p$ to be the set of
isomorphism classes of triples $(A,i,Z)$ over $\Fpb$ such
that $A\cong E\times E$ for some supersingular elliptic
curve $E$.
\item If $p\mid N^-$ define $\T_p$ to be the set of
isomorphism classes of triples $(A,i,Z)$ over $\Fpb$ such
that $A\cong E\times E$ for some supersingular elliptic curve
$E$ and $i$ is of type 1 or 2.
\item If $p\mid N^+$ define $\T_p$ to be the empty set.
\end{enumerate}
\end{definition}

\begin{remark}
In Proposition~\ref{crossing} we
will show that for $p\nmid N^+$ the set $\T_p$
parametrizes isomorphism classes of Eichler orders of type
$(N^+,pN^-)$.  Thus for each prime $p$ the set $\T_p$ is finite.
\medskip
\end{remark}

     Let $p\nmid N^+$ and let $(A,i,Z)\in\T_p$.
We now use the orientation $\{\phi_l\}_{l\mid N}$ on
$\ODelN$ to construct an orientation $\{\psi_l\}_{l\mid pN}$
on $\End(A,i,Z)$.
For $l\not=p$ let $T_l(A)$ be the $l$-adic Tate module of $A$
and let $U_l\supset T_l(A)$ be the lattice which corresponds
to the $l$-primary part of $Z$.  Then $U_l$ is free of rank 1
as a left module over $\ODelN\otimes\Zl$.  Using
Lemma~\ref{known} we get an isomorphism
\begin{equation} \label{isom}
\End(A,i,Z)\otimes\Zl\cong\End(U_l,i\otimes\Zl)\cong
\ODelN^{op}\otimes\Zl.
\end{equation}
Thus the orientation $\phi_l$ on $\ODelN\otimes\Zl$ induces
an orientation $\psi_l$ on $\End(A,i,Z)\otimes\Zl$.

     It remains to construct an orientation on
$\End(A,i,Z)\otimes\Zp$.  It follows from the computations
above (cf.\ (\ref{O1p}) and (\ref{cent})) that $\End(A,i,Z)$
acts on $\Lie(A)$ through scalar multiplication by elements
of $\Fps$.  If $p\nmid N$ this action gives an orientation
\begin{equation}
\psi_p:\End(A,i,Z)\otimes\Z_p\lra\Fps\cong R_{1,p}.
\end{equation}
If $p\mid N^-$ we note that $\Zps/p\Zps\cong\Fps$ embeds in
$\ODel\otimes\Fp$.  Let $e\in\Lie(A)$ be an eigenvector for
$\Fps$ such that
$(\ODel\otimes\Fp)\cdot v$ spans $\Lie(A)$.  (It follows from
(\ref{zps}) and (\ref{cases}) that $v$
exists and is uniquely determined up to scalar multiplication.)
For every $\alpha\in\End(A,i,Z)$ there is a unique
$\beta\in\ODel\otimes\Fp$ such that $\alpha(v)=i(\beta)\cdot v$.
Define $\overline{\psi}_p:\End(A,i,Z)\ra R_{1,p}$
by setting $\overline{\psi}_p(\alpha)=\phi_p(\beta)$.
Since $\End(A,i,Z)\otimes\Zp$ is an Eichler order of type
$(1,p^2)$, the homomorphism $\overline{\psi}_p$ lifts uniquely
to an orientation
\begin{equation}
\psi_p:\End(A,i,Z)\otimes\Zp\lra R_{1,p^2}.
\end{equation}

     Combining the above results we get the following
proposition:

\begin{prop} \label{orders}
(a) Let $t$ be a point of characteristic $p$ in the support of
the intersection of $\P_{D_1}$ with $\P_{D_2}$.
Then $(A_t,i_t,Z_t)\in\T_p$. \\[\smallskipamount]
(b) Let $(A,i,Z)\in\T_p$ and let $\{\psi_l\}_{l\mid pN}$
be the orientation on $\End(A,i,Z)$ induced by the
orientation $\{\phi_l\}_{l\mid N}$ on $\ODelN$.  Then
$(\End(A,i,Z),\{\psi_l\}_{l\mid pN})$ is an oriented Eichler
order of type $(N^+,pN^-)$.
\end{prop}

\section{Families of Eichler orders}
\label{families}
\setcounter{equation}{0}

     In this section we describe the relation between the
endomorphism rings of elements of $\T_p$ and isomorphism
classes of oriented Eichler orders of type $(N^+,pN^-)$.
It follows from Proposition~\ref{ss} that if $(A,i,Z)\in\T_p$
then $A\cong E\times E$ for any supersingular elliptic curve
$E$ over $\Fpb$.  As in \S\ref{points} we assume that $E$
is defined over $\Fp$.  We also let $F\in\Lambda=\End(E)$
denote the Frobenius endomorphism of $E$.

\begin{prop} \label{crossing}
(a) If $p\nmid N$ then the map $(A,i,Z)
\mapsto(\End(A,i,Z),\{\psi_l\}_{l\mid pN})$ gives a
bijection between $\T_p$ and the set of
isomorphism classes of oriented Eichler
orders of type $(N^+,pN^-)$. \\[\smallskipamount]
(b) If $p\mid N^-$ then the map $(A,i,Z)\mapsto
(\End(A,i,Z),\{\psi_l\}_{l\mid pN})$ gives a bijection
between the set of $(A,i,Z)\in\T_p$
such that $i$ has type 1,
and the set of isomorphism classes of oriented Eichler orders
of type $(N^+,pN^-)$.  The same statement holds for
$i$ of type 2.
\end{prop}

\proof We first show that the given maps are onto.  Let
$(\E,\{\mu_l\}_{l\mid pN})$ be an oriented Eichler order
of type $(N^+,pN^-)$.  Let $i_0:\ODel\ra\M_2(H_p)$ be an 
embedding.  Since the commutant of $i_0(\ODel)$ in $\M_2(H_p)$
is isomorphic to $\Delta(p)$, there is an embedding
$h_0:\E\ra\M_2(H_p)$ such that $h_0(\E)$ commutes with
$i_0(\ODel)$.  We will show that there exists $g\in\GL_2(H_p)$
and $Z\subset E\times E$
such that $gi_0(x)g^{-1}$ is a special embedding of $\ODel$
into $\M_2(\Lambda)\cong\End(E\times E)$ and
$gh_0(y)g^{-1}$ gives an isomorphism between the
oriented orders $(\E,\{\mu_l\}_{l\mid pN})$ and
$\End(E\times E,gi_0(x)g^{-1},Z)$.
Let $l$ be a prime such that $l\mid N$
and $l\not=p$.  Then $g^{-1}(T_l(E\times E))$ must satisfy
the following conditions:
\begin{enumerate}
\item The lattice $g^{-1}(T_l(E\times E))$
is stabilized by $i_0(\ODel)$ and $h_0(\E)$.
\item The orientations $\phi_l$ and $\mu_l$ are
compatible with respect to $g^{-1}(T_l(E\times E))$.
\end{enumerate}
To identify an appropriate $g$ we will first find lattices
which satisfy these properties.

     Suppose $l\mid N^+$.  Then
$\M_2(H_p)\otimes\Q_l\cong\M_4(\Q_l)$ acts on $\Q_l^4$,
and hence $i_0$ and $h_0$ induce an action of
$(\E\otimes\Z_l)\otimes(\ODel\otimes\Z_l)$ on $\Q_l^4$.
Using the Skolem-Noether theorem we may identify
$\E\otimes\Z_l$ with $\hat{\O}_{l^{n_l^+},1}$,
$\ODel\otimes\Z_l$ with $\M_2(\Z_l)$, and $\Q_l^4$ with
$\M_2(\Q_l)$ in such a way that this action
is isomorphic to $(A\otimes B)\cdot X=AXB^{\iota}$, where
$A\in\hat{\O}_{l^{n_l^+},1}$, $B\in\M_2(\Z_l)$, 
$X\in\M_2(\Q_l)$, and $\iota$ is the canonical involution of
$\M_2(\Q_l)$.  It follows that every lattice in
$\Q_l^4\cong\M_2(\Q_l)$
which is stabilized by $(\E\otimes\Z_l)\otimes(\ODel\otimes\Z_l)$
is a $\Q_l^{\times}$-multiple of one of the lattices
\begin{equation}
L_j=\left\{\begin{bmatrix}a&b\\c&d\end{bmatrix}\in\M_2(\Z_l):
l^j\mid c\text{ and }l^j\mid d\right\}
\end{equation}
for $0\le j\le n_l^+$.

     We may assume further that the identification of
$\ODel\otimes\Z_l$ with $\M_2(\Z_l)$ maps $\ODelN\otimes\Z_l$
onto $\hat{\O}_{l^{n_l^+},1}$.
Let $M_j\supset L_j$ be a lattice such that $M_j/L_j$ is
cyclic of order $l^{n_l^+}$ and $M_j$ is stabilized by both
$\E\otimes\Z_l$ and $\ODelN\otimes\Z_l$.  (The lattice $M_j$
corresponds to the $l$-primary part of $Z$.) Such a lattice
$M_j$ exists if and only if $j=0$ or $j=n_l^+$.  For those
cases we have
\begin{equation}
M_0=\left\{\begin{bmatrix}a&l^{-n_l^+}b\\c&d\end{bmatrix}:
a,b,c,d\in\Z_l\right\}
\end{equation}
and $M_{n_l^+}=\hat{\O}_{l^{n_l^+},1}$.  Note that $\phi_l$
induces opposite orientations on $\E\otimes\Z_l$ with respect
to the lattices $M_0$ and $M_{n_l^+}$ (cf.~(\ref{isom})).
Let $\m_l\subset\M_2(H_p)\otimes\Q_l$ denote the
stabilizer of the lattice $L_j$ ($j=0,n_l^+$) such that
$\phi_l$ induces the orientation $\mu_l$ on $\E\otimes\Z_l$
with respect to $M_j$.

     Suppose $l\mid N^-$ with $l\not=p$.  Then
$\M_2(H_p)\otimes\Q_l\cong\M_4(\Q_l)$ acts on
$\Q_l^4$, and hence $i_0$ and $h_0$ induce an action of
$(\E\otimes\Z_l)\otimes(\ODel\otimes\Z_l)$ on $\Q_l^4$.  Since
$\E\otimes\Z_l\cong\ODel\otimes\Z_l\cong\hat{\O}_{1,l}$, by
the Skolem-Noether theorem this action is isomorphic to
$(a\otimes b)\cdot x=axb^{\iota}$, where
$a,b\in\hat{\O}_{1,l}$, $x\in\hat{B}_l$, and $\iota$ is
the canonical involution of $\hat{B}_l$.  There are two
$\Q_l^{\times}$-equivalence classes of lattices in $\hat{B}_l$
stabilized by $(\E\otimes\Z_l)\otimes(\ODel\otimes\Z_l)$.
These are represented by $\hat{\O}_{1,l}$ and
$\pi\hat{\O}_{1,l}$, where $\pi$ is a uniformizer for
$\hat{\O}_{1,l}$.  The orientation $\phi_l$ on
$\ODelN\otimes\Z_l$ induces opposite orientations
on $\E\otimes\Z_l$ with respect to $\hat{\O}_{1,l}$ and
$\pi\hat{\O}_{1,l}$.  Choose $j=0,1$ so that $\phi_l$ induces
the orientation $\mu_l$ on $\E\otimes\Z_l$ with respect to
$\pi^j\hat{\O}_{1,l}$, and let $\m_l\subset\M_2(H_p)\otimes\Q_l$
denote the stabilizer of $\pi^j\hat{\O}_{1,l}$.

     Suppose $p\nmid N$.  Then by the Skolem-Noether theorem
we may identify
$\M_2(H_p)\otimes\Q_p$ with $\M_2(\hat{B}_p)$ in such a way
that $i_0(\ODel)\otimes\Z_p=\M_2(\Z_p)$ and
$h(\E)\otimes\Z_p=\hat{\O}_{1,p}\cdot I_2$. 
Suppose $p\mid N^-$ and $i$ has type 1.  Then we may identify
$\M_2(H_p)\otimes\Q_p$ with $\M_2(\hat{B}_p)$ so that
\begin{align} \label{ODel}
i_0(\ODel)\otimes\Z_p&=\left\{\begin{bmatrix}x&y\\py'&x'
\end{bmatrix}:x,y\in\Zps\right\} \\[.2cm]
h_0(\E)\otimes\Z_p&=\left\{\begin{bmatrix}a&b\pi\\pb\pi&a
\end{bmatrix}:a,b\in\Zps\right\} \label{E}
\end{align}
and $\phi_p$ induces the orientation $\psi_p$ on
$\E\otimes\Z_p$.  In either case we define $\m_p$ to be the
subring of $\M_2(H_p)\otimes\Q_p$ which corresponds to
$\M_2(\hat{\O}_{1,p})$ under this identification.

     Let $\m$ be a maximal order in $\M_2(H_p)$ such that
$\m\supset i_0(\ODel)$, $\m\supset h_0(\E)$, and
$\m\otimes\Z_l=\m_l$ for all
$l$ such that $l\mid pN$.  Let $\c$ denote the commutant of
$i_0(\ODel)$ in $\m$.  It follows from the constructions above
that if $l\mid pN$ then $\c\otimes\Z_l=h_0(\E)\otimes\Z_l$.
In particular, if $l\mid N$ and $l\not=p$ then
$\c\otimes\Z_l$ is a
local Eichler order of type $(1,l^{n_l^-})$; if $p\nmid N$
then $\c\otimes\Z_p$ is a local Eichler order of type $(1,p)$;
and if $p\mid N^-$ then $\c\otimes\Z_p$ is a local Eichler
order of type $(1,p^2)$. 
Suppose $l\nmid pN$.  Then $\m\otimes\Z_l\cong\M_4(\Z_l)$ and
$\ODel\otimes\Z_l\cong\M_2(\Z_l)$, and hence
$\c\otimes\Z_l\cong\M_2(\Z_l)$.  Combining these local
results we conclude that $\c$ is an Eichler order of type
$(1,pN^-)$.

     It follows from Theorems~21.6 and 34.9 in \cite{mo}
that the maximal orders of $\M_2(H_p)$ are all conjugate.
Therefore there exists $g\in\GL_2(H_p)$ such that
$g\m g^{-1}=\M_2(\Lambda)$.  Let $i(x)=gi_0(x)g^{-1}$ and
$h(y)=gh_0(y)g^{-1}$ for $x\in\ODel$ and $y\in\E$.  By the
construction of $\m$ we see that $i(\ODel)$ and $h(\E)$ are
contained
in $\M_2(\Lambda)\cong\End(E\times E)$.  In addition, if
$p\mid N^-$ then it follows from (\ref{ODel}) that
$i:\ODel\ra\End(E\times E)$ is a special embedding of type 1.
Let $Z$ be the unique cyclic subgroup of $E\times E$ of order
$N^+$ which is stabilized by both $h(\E)$ and $i(\ODelN)$.
Let $\{\psi_l\}_{l\mid pN}$ be the orientation on
$\E\cong\End(E\times E,i,Z)$ induced by $\{\phi_l\}_{l\mid N}$.
Then by the constructions above we have $\psi_l=\mu_l$ for
all $l$ such that $l\mid N$.  Thus if $p\mid N$ then $h$
induces an isomorphism
\begin{equation} \label{orisom}
(\E,\{\mu_l\}_{l\mid pN})\cong
(\End(E\times E,i,Z),\{\psi_l\}_{l\mid pN})
\end{equation}
as required.
Suppose $p\nmid N$.  If $\psi_p=\mu_p$ then (\ref{orisom})
still holds, while if $\psi_p\not=\mu_p$ then (\ref{orisom})
holds after we replace $i(x)$ with $Fi(x)F^{-1}$ and $h(x)$
with $Fh(x)F^{-1}$.

     To prove that our map is one-to-one we need to show that
if the endomorphism rings of the triples $(A_1,i_1,Z_1)$ and
$(A_2,i_2,Z_2)$ are isomorphic as oriented orders then
$(A_1,i_1,Z_1)\cong(A_2,i_2,Z_2)$.  We may assume that
$A_1=A_2=E\times E$.  Then there is an oriented Eichler
order $(\E,\{\psi_l\}_{l\mid pN})$ in $\Delta(p)$ of type
$(N^+,pN^-)$, and embeddings $h_j:\E\ra\M_2(\Lambda)$ for
$j=1,2$ which induce isomorphisms between the oriented orders
$(\E,\{\psi_l\}_{l\mid pN})$ and $\End(E\times E,i_j,Z_j)$.
We need to show there is $g\in\GL_2(\Lambda)$ such that
$Z_2=gZ_1$ and $i_2(x)=gi_1(x)g^{-1}$ for all $x\in\ODel$.

     By the Skolem-Noether theorem there is $g\in \GL_2(H_p)$
such that $i_2(x)=gi_1(x)g^{-1}$ and $h_2(y)=gh_1(y)g^{-1}$
for all $x\in\ODel$ and $y\in\E$.  We may assume that
$g\in\M_2(\Lambda)$, and that $|\Nr(g)|$ is as small as
possible.  We claim that $g\in\GL_2(\Lambda)$.
Let $R$ be the order in $\M_2(\Lambda)$ generated by
$i_2(\ODel)$ and $h_2(\E)$.  The lattices
$\Lambda^2$ and $g(\Lambda^2)$ are stabilized by
$i_2(\ODel)=gi_1(\ODel)g^{-1}$ and $h_2(\E)=gh_1(\E)g^{-1}$,
and hence also by $R$.  Let $l\not=p$.  Then $\ODel\otimes\Z_l$,
$\E\otimes\Z_l$, and $R\otimes\Z_l$ stabilize
$T_l=T_l(E\times E)$ and $g(T_l)$.  Suppose $l\mid N$
and $l\not=p$.  Then there is a lattice $U_l\supset T_l$ such
that $U_l/T_l$ is cyclic of order $l^{n_l^+}$ and both $U_l$
and $g(U_l)$ are stabilized by $i_2(\ODelN)\otimes\Z_l$
and $h_2(\E)\otimes\Z_l$.  The orientations on
$\E\otimes\Z_l\cong h_2(\E)\otimes\Z_l$ induced by the
orientation $\phi_l$ on $\ODelN\otimes\Z_l\cong
i_2(\ODelN)\otimes\Z_l$ with respect to the
lattices $U_l$ and $g(U_l)$ are identical.

     Suppose $l\mid N$ with $l\not=p$.  We saw above that there
are two $\Q_l^{\times}$-equivalence classes of pairs of lattices
$M_l\supset L_l$ in $\Q_l^4$ such that $M_l/L_l$ is cyclic of
order $l^{n_l^+}$, $M_l$ is stabilized by
$\ODelN\otimes\Z_l$ and $\E\otimes\Z_l$, and $L_l$ is
stabilized by $\ODel\otimes\Z_l$ and $\E\otimes\Z_l$.
The orientation $\phi_l$ induces opposite orientations on
$\E\otimes\Z_l$ with respect to the two possibilities for $M_l$.
Therefore $g(T_l)$ lies in the same class as $T_l$, so
$g(T_l)=l^kT_l$ for some
$k\ge0$.  Suppose $l\nmid pN$.  Then $\ODel\otimes\Z_l\cong
\E\otimes\Z_l\cong\M_2(\Z_l)$ and hence
$R\otimes\Z_l\cong\M_4(\Z_l)$ stabilizes $T_l$ and
$g(T_l)$.  It follows that $g(T_l)=l^kT_l$ for some $k\ge0$.
In both cases we have $g=l^k\cdot g_0$ for
some $g_0\in\M_2(\Lambda)$.  By the
minimality of $|\Nr(g)|$ we get $k=0$, and hence
$g(T_l)=T_l$.  It follows that
$g\in\GL_2(\Lambda\otimes\Z_l)$ for every $l\not=p$.

     Suppose $p\nmid N$.  Then
$\ODel\otimes\Z_p\cong\M_2(\Z_p)$ and
$\E\otimes\Z_p\cong\hat{\O}_{1,p}$, and hence
$R\otimes\Z_p\cong\M_2(\hat{\O}_{1,p})$.  Therefore
$g(\Lambda^2)\otimes\Z_p$ is stabilized by
$\M_2(\Lambda\otimes\Z_p)$, so
$g(\Lambda^2)\otimes\Z_p=F^k\cdot(\Lambda^2\otimes\Z_p)$
for some $k\ge0$.  It follows from the minimality of $|\Nr(g)|$
that $k=0,1$.  However, if
$g(\Lambda^2)\otimes\Z_p=F\cdot(\Lambda^2\otimes\Z_p)$ then
$g\in F\cdot\GL_2(\Lambda)$, and hence $i_1(x)$ and
$i_2(x)=gi_1(x)g^{-1}$ can't both be of type 1,
contrary to assumption.  Hence
$g\in\GL_2(\Lambda\otimes\Z_p)$.

     Suppose $p\mid N^-$.  It follows using (\ref{ODel}) and
(\ref{E}) that we may identify $\M_2(\Lambda\otimes\Z_p)$
with $\M_2(\hat{\O}_{1,p})$ in such a way that
$R\otimes\Z_p$ contains the
matrices $\begin{bmatrix}1&0\\0&0\end{bmatrix}$ and
$\begin{bmatrix}0&1\\p&0\end{bmatrix}$.  Since
$g(\Lambda^2)\otimes\Z_p$ is an $\hat{\O}_{1,p}$-lattice in
$(H_p\otimes\Qp)^2\cong\hat{B}_p^2$ which is stabilized by
$R\otimes\Z_p$, it is $\Qp^{\times}$-equivalent to one of the
lattices $\hat{\O}_{1,p}\oplus\hat{\O}_{1,p}$,
$\hat{\O}_{1,p}\oplus\pi\hat{\O}_{1,p}$,
$\hat{\O}_{1,p}\oplus p\hat{\O}_{1,p}$,
$\pi\hat{\O}_{1,p}\oplus\pi\hat{\O}_{1,p}$,
$\pi\hat{\O}_{1,p}\oplus p\hat{\O}_{1,p}$, or
$\pi\hat{\O}_{1,p}\oplus p\pi\hat{\O}_{1,p}$.
However, since $i_1$ and $i_2$ are both special embeddings,
$g(\Lambda^2)\otimes\Z_p$ is not $\Qp^{\times}$-equivalent
to either $\hat{\O}_{1,p}\oplus\pi\hat{\O}_{1,p}$ or
$\pi\hat{\O}_{1,p}\oplus p\hat{\O}_{1,p}$.  Since $i_1$ and
$i_2$ are both of of type 1, $g(\Lambda^2)\otimes\Z_p$ is
not equivalent to either
$\hat{\O}_{1,p}\oplus p\hat{\O}_{1,p}$ or
$\pi\hat{\O}_{1,p}\oplus\pi\hat{\O}_{1,p}$.  Since
$\phi_p$ induces the same orientation on $\E\otimes\Z_p$
with respect to the embedding pairs $(i_1,h_1)$ and
$(i_2,h_2)$, $g(\Lambda^2)\otimes\Z_p$ is not equivalent to
$\pi\hat{\O}_{1,p}\oplus p\pi\hat{\O}_{1,p}$.  We conclude
that $g(\Lambda^2)$ is $\Qp^{\times}$-equivalent to
$\hat{\O}_{1,p}\oplus\hat{\O}_{1,p}$.  Hence
$g(\Lambda^2)\otimes\Z_p=p^k\cdot(\Lambda^2\otimes\Z_p)$ for
some $k\ge0$.  By the minimality of $|\Nr(g)|$ we get $k=0$,
and hence $g\in\GL_2(\Lambda\otimes\Z_p)$.

     Combining the above results we get $g\in\GL_2(\Lambda)$.
Finally, the $l$-primary subgroup of $Z_1$ is the unique cyclic
subgroup of $E\times E$ of order $l^{n_l^+}$ which is
stabilized by both $h_1(\E)$ and $i_1(\ODelN)$.  A
similar statement holds for $j=2$, and
so $g$ maps the $l$-primary subgroup of
$Z_1$ to the $l$-primary subgroup of
$Z_2$.  Therefore $Z_2=gZ_1$.~\qed

     Recall that $\G_p$ is a set of representatives of the
proper equivalence classes of a genus of quaternary quadratic
spaces over $\Z$, and that $\s_p$ denotes the set of Eichler
orders in $\Delta(p)$ of type $(N^+,pN^-)$.
We may view the elements of $\s_p$ as quadratic spaces over
$\Z$ with the $\Z$-valued quadratic
forms induced by the reduced norm on $\Delta(p)$.

\begin{prop} \label{S}
(a) Every $L\in\G_p$ which represents 1
over $\Z$ is properly equivalent to some $\E\in\s_p$.
\\[\smallskipamount]
(b) The orders $\E,\E'\in\s_p$ are properly equivalent
quadratic spaces if
and only if $\E'=a\E a^{-1}$ for some $a\in\Delta(p)^{\times}$.
\end{prop}

     To facilitate the proof we recall the following
well-known fact (cf.~\cite[I,\,Th.\,3.3]{vig}):

\begin{lemma} \label{self}
Let $B$ be a quaternion algebra over a field $F$ whose
characteristic is not 2.  Then
every proper self-isometry of $B$ has the form
$\phi(x)=axb^{-1}$ for some $a,b\in B^{\times}$ such that
$\Nr(a)=\Nr(b)$.
\end{lemma}

\noindent{\em Proof of Proposition \ref{S}:} (a) Choose an
Eichler order $\E_0\subset\Delta(p)$ of type $(N^+,pN^-)$.
Then $\E_0$ with the reduced norm form belongs to the same
genus as $L$, so by the weak Hasse principle
\cite[p.\,76]{ca} $L$ is isometric to a lattice $J$
in $\Delta(p)$.  Let $l$ be a prime.  Then $\E_0\otimes\Z_l$
is isomorphic to one of the standard local Eichler orders
given in (\ref{standard}) and (\ref{standard2}).
Since the standard local Eichler orders are
stabilized by the appropriate canonical involutions,
there is a {\em proper} isometry between
$J\otimes\Z_l$ to $\E_0\otimes\Z_l$.  It follows from
Lemma~\ref{self} that there are
$a_l,b_l\in(\Delta(p)\otimes\Q_l)^{\times}$ such that
$J\otimes\Z_l=a_l(\E_0\otimes\Z_l)b_l^{-1}$.

     Let $\E_J\subset\Delta(p)$ be the right order of $J$.
Since $\E_J\otimes\Z_l=b_l(\E_0\otimes\Z_l)b_l^{-1}$ for
every prime $l$ we have $\E_J\in\s_p$.  Let $u\in J$ be such
that $\Nr(u)=1$.  Then $u\E_J\subset J$, and for every prime
$l$ the quadratic spaces
$u\E_J\otimes\Z_l$ and $J\otimes\Z_l$ are both isometric to
the local Eichler order $\E_0\otimes\Z_l$.  Since
$\E_0\otimes\Z_p$ is a nondegenerate quadratic space, this
implies $u\E_J\otimes\Z_l=J\otimes\Z_l$.  Thus $u\E_J=J$,
and hence $J$ and $L$ are isometric to $\E_J\in\s_p$.
If the isometry $f:L\ra\E_J$ is not proper, replace $f$ with
$\iota\circ f$ and $\E_J$ with $\iota(\E_J)$, where $\iota$
is the canonical involution of $\Delta(p)$. \\[\smallskipamount]
(b) If $\E'=a\E a^{-1}$ with
$a\in\Delta(p)^{\times}$ then $x\mapsto
axa^{-1}$ is a proper isometry from $\E$ to
$\E'$.  Conversely,
suppose $\phi:\E\ra \E'$ is a proper isometry.
By Lemma~\ref{self} we have 
$\phi(x)=axb^{-1}$ for some $a,b\in\Delta(p)^{\times}$.
Clearly $\phi(1)=ab^{-1}$ is a unit in $\E'$, so
$a\E a^{-1}=\phi(\E)(ab^{-1})^{-1}=\E'$.
\Qed

\section{Universal deformations}
\label{universal} \setcounter{equation}{0}

     Let $(A_t,i_t,Z_t,\beta_t)$ correspond to a point
$t\in\X(\Fpb)$.  Let $W_p=W(\Fpb)$ denote the ring of Witt
vectors with coefficients in $\Fpb$, and let $\hat{\X}_t$
be the completion of $\X\otimes W_p$ at $t$.  Since $\X$
is a fine moduli space, $\hat{\X}_t$ is a universal
deformation space for $(A_t,i_t,Z_t,\beta_t)$.
In this section we study $\hat{\X}_t$.
Let $\hat{A}_t$ be the formal group
of $A_t$ and let $\hat{\imath}_t:\ODel\otimes\Zp \ra\End(\hat{A}_t)$
be the map induced by $i_t$.  To begin we show that $\hat{\X}_t$
is a universal deformation space for the pair
$(\hat{A}_t,\hat{\imath}_t)$.

\begin{lemma} \label{equiv}
Let $p$ be a prime such that $p\nmid mN^+$, and let
$R$ be a complete Noetherian local ring with residue field
$\Fpb$.  Then there are natural bijections between \\[.2cm]
(a) The set of deformations over $R$ of
the 4-tuple $(A_t,i_t,Z_t,\beta_t)$, \\[\smallskipamount]
(b) The set of deformations over $R$ of
the pair $(A_t,i_t)$, and \\[\smallskipamount]
(c) The set of deformations over $R$ of
the pair $(\hat{A}_t,\hat{\imath}_t)$.
\\[.2cm]
Therefore $\hat{\X}_t$ serves as a universal
deformation space for either $(A_t,i_t,Z_t,\beta_t)$,
$(A_t,i_t)$, or $(\hat{A}_t,\hat{\imath}_t)$.
\end{lemma}

\proof Since $p\nmid mN^+$, Hensel's Lemma
\cite[I,\,\,Th.\,4.2(d)]{mi} implies that if $(A,i)$
is a deformation of $(A_t,i_t)$ defined over
$R$ then $Z_t$ extends uniquely to a subgroup scheme
$Z\subset A$ which is finite and flat of order $N^+$ over
$R$, and $\beta_t$ extends uniquely to a $\Gamma_1(m)$-structure
$\beta$ on $A$.  This gives a bijection between (a) and (b).
The Serre-Tate lifting theorem (see the appendix to
\cite{cd}) gives a bijection between (b) and (c).~\Qed

\begin{prop} \label{univdef}
Let $p$ be a prime such that $p\nmid mN^+$, and let $t\in\X(\Fpb)$.
\\[.2cm]
(a) If $p\nmid N^-$ then $\hat{\X}_t\cong\Spf\:W[[u]]$. \\[\smallskipamount]
(b) If $p\mid N^-$ and $i$ is of type 1 or 2
then $\hat{\X}_t\cong\Spf\:W[[u]]$.  If $p\mid
N^-$ and $i$ is of type 3 then
$\hat{\X}_t\cong\Spf\:W[[u,v]]/(uv-p)$.
\end{prop}

\proof (a) In this case $\hat{A}_t$ has multiplication by
$\ODel\otimes\Zp\cong\M_2(\Zp)$ and must
therefore be isomorphic to a product
$G_0\times G_0$, where $G_0$ is a formal
group of dimension 1 and height 2 over
$\Fpb$.  By the same reasoning any
deformation of $\hat{A}_t$ with
multiplication by $\ODel\otimes\Zp$ has
the form $G\times G$, where $G$ is a
deformation of $G_0$, and conversely any
deformation $G$ of $G_0$ gives a unique
deformation $G\times G$ of $\hat{A}_t$
with multiplication by $\ODel\otimes\Zp$.
Therefore deformations of the pair
$(\hat{A}_t,\hat{\imath}_t)$ are
equivalent to deformations of $G_0$.  In \cite{lt} it is
proved that the universal deformation space of the
formal group $G_0$ is $\Spf\:W_p[[u]]$.  Hence by
Lemma~\ref{equiv} we get $\hat{\X}_t\cong\Spf\:W_p[[u]]$. \\[\smallskipamount]
(b) In this case $\hat{A}_t$ has multiplication by
$\ODel\otimes\Zp\cong\hat{\O}_{1,p}$.
We use Drinfeld's theory in \cite{cd} to
interpret the formal scheme $\hat{\H}_p$
(the ``$p$-adic upper half-plane'')
as a moduli space for rigidified formal
groups of dimension 2 and height 4 with
a special action by $\hat{\O}_{1,p}$.  Associated to the pair
$(\hat{A}_t,\hat{\imath}_t)$ is an
equivalence class of closed points in $\hat{\H}_p$.  The
formal neighborhood of any of these
points is a universal deformation space
for $(\hat{A}_t,\hat{\imath}_t)$, and hence also for
$(A_t,i_t,Z_t,\beta_t)$.

     To determine the structure of this formal neighborhood
we use the explicit description of $\hat{\H}_p$
found in \cite[pp.\,650--51]{jt} and
\cite[I,\,\,\S3]{bc}.  The special fiber
of $\hat{\H}_p$ is an infinite tree
consisting of projective lines which meet transversely at
their $\Fp$-rational points.  The formal
neighborhood of a closed point in $\hat{\H}_p$
takes two different forms, depending on
whether or not it is a crossing point of
the special fiber of $\hat{\H}_p$.
To distinguish the crossing points from
the other points in the special fiber of
$\hat{\H}$ we let $\pi$ be an element of
$\ODel\otimes\Zp$ such that $\pi^2=p$.
Then $\pi$ is a generator for the
maximal ideal of $\ODel\otimes\Zp$,
and $\pi$ induces an endomorphism of
$\Lie(\hat{A}_t)$ whose square is zero.  If $i_t$ is
of type 1 or 2 then $\pi$ induces a non-trivial
endomorphism of $\Lie(\hat{A}_t)$.  In this case $t$
is not a crossing point of the special
fiber of $\hat{\H}_p$, and \cite[p.\,650]{jt} gives
$\hat{\X}_t\cong\Spf\:W_p[[u]]$.  If $i_t$ is of
type 3 then $\pi$ induces the zero map on $\Lie(\hat{A}_t)$.
In this case $t$ is a crossing point of the
special fiber of $\hat{\H}_p$, and by \cite[p.\,650]{jt} we
have $\hat{\X}_t\cong\Spf\:W_p[[u,v]]/(uv-p)$.~\qed \medskip

     The following consequence of Proposition~\ref{univdef}
is presumably well-known:

\begin{cor} \label{regular}
$\X\otimes\Z[1/mN^+]$ is a regular scheme.
\end{cor}

     Let $R$ be a complete noetherian local
ring with residue field $\Fpb$ and let
$(\hat{A},\hat{\imath})$ be a deformation
of $(\hat{A}_t,\hat{\imath}_t)$ defined over
$R$.  The reduction map $R\ra\Fpb$ induces inclusions
$\End(\hat{A})\subset\End(\hat{A}_t)$ and
$\End(\hat{A},\hat{\imath})\subset\End(\hat{A}_t,\hat{\imath}_t)$.
Let $\univ$ be a universal
deformation of $(\hat{A}_t,\hat{\imath}_t)$
defined over $\hat{\X}_t$, and let $S$ be a
subset of $\End(\hat{A}_t,\hat{\imath}_t)$.
We define $\hat{\X}_t(S)$ to be the largest formal subscheme of
$\hat{\X}_t$ such that $S\subset\End(\underline{\hat{A}}
\times_{\hat{\X}_t}\hat{\X}_t(S))$.
The following facts are easily verified:

\begin{lemma} \label{facts}
(a) $\hat{\X}_t(S)$ is closed in $\hat{\X}_t$.
\\[\smallskipamount]
(b) $\hat{\X}_t(S_1\cup S_2)=\hat{\X}_t(S_1)
\cap\hat{\X}_t(S_2)$. \\[\smallskipamount]
(c) Let $\Zp[S]$ be the $\Zp$-subalgebra
of $\End(\hat{A}_t,\hat{\imath}_t)$ generated by $S$.  If
$S\subset S'\subset\Zp[S]$ then $\hat{\X}_t(S')=\hat{\X}_t(S)$.
\end{lemma}

     Let $\gamma_1,\gamma_2\in \End(\hat{A}_t,\hat{\imath}_t)$.
The intersection multiplicity
$\alpha_p(\gamma_1,\gamma_2)$ of $\hat{\X}_t(\gamma_1)$
with $\hat{\X}_t(\gamma_2)$ is defined to be the
$W_p$-length of the coordinate ring of
$\hat{\X}_t(\gamma_1)\cap\hat{\X}_t(\gamma_2)=
\hat{\X}_t(\{\gamma_1,\gamma_2\})$.
To compute the arithmetic intersection numbers of our
divisors we need to evaluate $\alpha_p(\gamma_1,\gamma_2)$
for certain $\gamma_1,\gamma_2\in\End(\hat{A}_t,\hat{\imath}_t)$.

    Assume first that $p$ is not ramified in $\Delta$.  Then
$\hat{\X}_t$ is a universal deformation space for a formal
group $G_0$ over $\Fpb$ of dimension 1 and height 2.  By
Proposition~\ref{orders} we have
$\End(\hat{A}_t,\hat{\imath}_t)\cong\End(G_0)\cong\hat{\O}_{1,p}$.
The intersection multiplicity $\alpha_p(\gamma_1,\gamma_2)$
may be computed using the formulas in
\cite[Prop.\,5.4]{gk}.  In order to
state these formulas we define a
quadratic form over $\Zp$,
\begin{equation}
Q(x,y,z)=\Nr(x+y\gamma_1+z\gamma_2),
\end{equation}
where $\Nr$ is the reduced norm form on
$\End(\hat{A}_t,\hat{\imath}_t)\cong\hat{\O}_{1,p}$.
We wish to define invariants $a_1,a_2,a_3$ of $Q$.
If $p>2$ we diagonalize $Q$ over $\Zp$ and define
$a_1\leq a_2\leq a_3$ to be the $p$-adic
valuations of the coefficients of the
diagonal form of $Q$.  If $p=2$
the definition of $a_1,a_2,a_3$ is more
complicated and may be found in
\cite[\S4]{gk}.  In either case we have $a_1=0$, so by
\cite[Prop.\,5.4]{gk} we get
\begin{equation} \label{unram}
\alpha_p(\gamma_1,\gamma_2)=
\begin{cases}
\dst
\frac{a_3-a_2+1}{2}p^{a_2/2}+
\sum_{i=0}^{(a_2-2)/2}(a_2+a_3-4i)p^i&
\mbox{ if }a_2\equiv0\;\pmod{2}, \\[.5cm]
\dst\sum_{i=0}^{(a_2-1)/2}
(a_2+a_3-4i)p^i&\mbox{ if }a_2\equiv1\;\pmod{2}.
\end{cases}
\end{equation}
Note that the intersection multiplicity
depends only on the $\Zp$-isometry
class of $Q$, and not on the
particular $\gamma_1,\gamma_2$ that were used to
define $Q$.  Therefore we may write
$\alpha_p(Q)=\alpha_p(\gamma_1,\gamma_2)$.

     Now suppose that $p\mid N^-$ is ramified in $\Delta$.
We may assume without loss of generality that $i_t$ is of
type 1.  By Proposition~\ref{orders} we see that
$\End(\hat{A}_t,\hat{\imath}_t)\cong\End(A_t,i_t,Z_t)\otimes\Zp$
is a local Eichler order of type $(1,p^2)$.  Suppose that
for $j=1,2$ there are embeddings of $\O_{D_j}$ into
$\End(A_t,i_t,Z_t)$ with image $\Z[\gamma_j]$.  Using
(\ref{assume2}) we may assume that $p\nmid D_1$.  It follows
then from Remark~\ref{split} that
$\O_{D_1}\otimes\Z_p\cong\Zp[\gamma_1]\cong\Zps$.
Let $(\hat{A},\hat{\imath})$ be a
deformation of $(\hat{A}_t,\hat{\imath}_t)$
such that $\gamma_1\in\End(\hat{A},\hat{\imath})$.  Since
$i_t$ is of type
1, $\End(\hat{A})$ contains a subring which is conjugate in
$\End(\hat{A}_t)\cong\M_2(\hat{\O}_{1,p})$ to
\begin{equation} \label{tensor}
\left\{\left[\begin{array}{cc}a&b\\pc&d\end{array}
\right]:a,b,c,d\in\Zps\right\}.
\end{equation}
Therefore we are in a situation much
like the case where $p$ is unramified in
$\Delta$: There is a formal group $G_0$
over $\Fpb$ of dimension 1 and height 2,
a map $\tau:\Zps\ra\End(G_0)$,
and a deformation $G$ of $G_0$ such
that $\hat{A}\cong G\times G$ and
$\tau(\Zps)\subset\End(G)$.
Conversely, any deformation $G$ of
$G_0$ such that $\tau(\Zps)\subset\End(G)$
gives a deformation
$(\hat{A},\hat{\imath})$ of
$(\hat{A}_t,\hat{\imath}_t)$ such that $\gamma_1\in
\End(\hat{A},\hat{\imath})$.
It follows that $\hat{\X}_t(\gamma_1)$ is
isomorphic to the universal deformation
space of the formal $\Zps$-module $(G_0,\tau)$.  It is
proved in \cite{qc} that the
universal deformation of $(G_0,\tau)$ is
the canonical lifting $\uG$ of $G_0$ associated
to $\tau$, which is defined over $W_p$.  Therefore we have
$\hat{\X}_t(\gamma_1)\cong\Spf\:W_p$ and 
$\hat{\X}_t(\{\gamma_1,\gamma_2\})\cong\Spf(W_p/(p^{k+1}))$
for some $k\geq0$.

     To determine $k$ we use an indirect argument.  Let
$(\hat{A},\hat{\imath})$ be the restriction of $\univ$
to $\hat{\X}_t(\gamma_1)\cong\Spf\:W_p$, and for
$n\ge0$ let $(\hat{A}_n,\hat{\imath}_n)$ be the restriction
of $\univ$ to $\Spf(W_p/(p^{n+1}))$.
Then $\hat{A}_n\cong G_n\times G_n$, where
$G_n=\uG\otimes(W_p/p^{n+1}W_p)$.  By
\cite[Prop.\,\,3.3]{qc} we have
\begin{equation} \label{endo}
\End(G_n)=\tau(\Zps)+p^n\hat{\O}_{1,p}.
\end{equation}
Using (\ref{cent}) we get
\begin{equation} \label{isomorphism}
\End(\hat{A}_n,\hat{\imath}_n)=\left\{\left[
\begin{array}{cc}
a&b\pi\\pb\pi&a\end{array}\right]:
a\in\tau(\Zps),\;b\in p^n\tau(\Zps)\right\}.
\end{equation}
It follows that $\End(\hat{A}_n,\hat{\imath}_n)$ is a local
Eichler order of type $(1,p^{2n+2})$.

     Recall that $k$ is the largest value of
$n$ such that $\Zp[\gamma_1,\gamma_2]$ is contained in
$\End(\hat{A}_n,\hat{\imath}_n)$.  Since $\Z_p[\gamma_1,\gamma_2]$
and $\End(\hat{A}_n,\hat{\imath}_n)$ are both local
Eichler orders which contain $\Zp[\gamma_1]\cong\Zps$, this
implies
$\Zp[\gamma_1,\gamma_2]=\End(\hat{A}_k,\hat{\imath}_k)$.
It follows that the reduced discriminant of
$\Zp[\gamma_1,\gamma_2]$ is $p^{2k+2}$.  Let $\delta$ denote
the reduced discriminant of the global order
$\Z[\gamma_1,\gamma_2]$.  Then we have $v_p(\delta)=2k+2$,
and hence
\begin{equation} \label{ram}
\alpha_p(\gamma_1,\gamma_2)=k+1=\frac{1}{2}\cdot v_p(\delta).
\end{equation}

     The following lemma will be used in determining the
relationship between $\hat{\X}_t(\sigma(\O_D))$ and the
divisor $\P_{D,\pm b}$.

\begin{lemma} \label{reduced}
Let $G_0$ be a formal group of dimension 1
and height 2 over $\Fpb$ and let $\hatU$
be a universal deformation space for $G_0$.
Let $R_p$ be an order in a quadratic extension $K$ of $\Qp$,
and let $\phi:R_p\ra\End(G_0)$ be a ring homomorphism.
Then $\hatU(\phi(R_p))$ is reduced.
\end{lemma}

\proof We have $R_p=\Zp+\Zp a$ for some
$a\in R_p$, so by Lemma~\ref{facts}(a) we
get $\hatU(\phi(a))=\hatU(\phi(R_p))$.
The subscheme $\hatU(\phi(a))$ of
$\hatU\cong\Spf\:W[[u]]$ is defined by an
equation of the form $f(u)-u=0$, where
$f(u)\in W_p[[u]]$ (cf.~\cite[p.\,58]{lt}).
Let $p^k$ be the conductor of the order $R_p$.
Then the Weierstrass degree of $f(u)-u$ is
computed in \cite[Th.\,1.1]{kk1} to be
\begin{equation} \begin{array}{rl}
p^k+2p^{k-1}+\cdots+2p+2&\mbox{if $K/\Qp$ is unramified,}
\\[.2cm]
2p^k+2p^{k-1}+\cdots+2p+2&\mbox{if $K/\Qp$ is ramified.}
\end{array} \end{equation}
The power series $f(u)-u$ is divisible by an
irreducible factor corresponding to a
quasi-canonical lifting of level $l$ for
each $0\leq l\leq k$.  The Weierstrass
degrees of these factors are computed in
\cite[Prop.\,5.3]{qc} to be
\begin{equation} \begin{array}{cl}
1&\mbox{if $K/\Qp$ is unramified and $l=0$,}
\\[.2cm]
p^l+p^{l-1}&\mbox{if $K/\Qp$ is unramified and $l\geq1$,}
\\[.2cm]
2p^l&\mbox{if $K/\Qp$ is ramified.}
\end{array} \end{equation}
Comparing Weierstrass degrees we find that $f(u)-u$ is the
product of the quasi-canonical lifting
factors for $0\leq l\leq k$ and a unit power series.  Since
the quasi-canonical lifting factors are
irreducible and have different degrees, the quotient
$W_p[[u]]/(f(u)-u)$ is reduced, as
claimed. \hspace{\fill}$\square$ \medskip

     Choose $D$ and $b=(b_l)_{l\mid N}$ as in \S\ref{heeg}.
Let $t\in\X(\Fpb)$ be such that $(A_t,i_t,Z_t)\in\T_p$.
We wish to consider the restriction $\hat{\P}_{D,\pm b}$
of the divisor $\P_{D,\pm b}$ to the completion
$\hat{\X}_t$ of $\X\otimes W_p$ at $t$.
Let $M_p$ denote the field of fractions of
$W_p$ and replace $X$ with $X\otimes M_p$ and
$P_{D,\pm b}$ with $P_{D,\pm b}\otimes M_p$.
Then $\P_{D,\pm b}\otimes W_p$ is the closure of
$P_{D,\pm b}\otimes M_p$ in $\X\otimes W_p$.  
For each $x$ in the support of $P_{D,\pm b}\otimes M_p$
define $\xhat$ to be the closure of $x$ in $\hat{\X}_t$.
Then we have $\hat{\P}_{D,\pm b}=\sum\:(\hat{x})$,
where the sum is taken over all points $x$ in the support
of $P_{D,\pm b}\otimes M_p$ such that $t\in\xbar$. 

     Let $\univZ$ be a universal deformation of
$(A_t,i_t,Z_t)$ defined over $\hat{\X}_t$, and let
$\{\psi_l\}_{l\mid pN}$ be the orientation on
$\End(A_t,i_t,Z_t)$ induced by $\{\phi_l\}_{l\mid N}$.
Define a $b$-embedding to be a ring homomorphism
$\sigma:\O_{D}\ra\End(A_t,i_t,Z_t)$ such that
$\psi_l\circ\sigma(\sqrt{D})=b_l$ for all $l\mid N$.
Say that $\sigma$ is a $\pm b$-embedding if $\sigma$ is
either a $b$-embedding or a $(-b)$-embedding.  The relation
between $\hat{\X}_t(\sigma(\O_{D}))$ and
$\hat{\P}_{D,\pm b}$ is given by the following lemma.

\begin{lemma} \label{lift}
(a) If $\sigma$ is a $\pm b$-embedding then
$\hat{\X}_t(\sigma(\O_{D}))$
is contained in the support of $\hat{\P}_{D,\pm b}$.
\\[\smallskipamount]
(b) Every irreducible component
of the support of $\hat{\P}_{D,\pm b}$ lies in
$\hat{\X}_t(\sigma(\O_{D}))$ for some $\pm b$-embedding
$\sigma$.
\\[\smallskipamount]
(c) Let $\sigma,\sigma'$ be $\pm b$-embeddings.  Then the
components in the support of
$\hat{\X}_t(\sigma(\O_{D}))$ and
$\hat{\X}_t(\sigma'(\O_{D}))$ are all
different unless $\sigma(\O_{D})=
\sigma'(\O_{D})$.
\end{lemma}

\proof (a) Suppose $p\nmid N^-$.  In the proof of
Proposition~\ref{univdef}\,(a) we
showed that $\hat{\X}_t$ is a universal
deformation space for a formal group
$G_0$ of dimension 1 and height 2.  The
map $\sigma$ induces an embedding of
$\O_{D}\otimes\Zp$ into $\End(G_0)$,
and $\hat{\X}_t(\sigma(\O_{D}))$ is
defined by the requirement that the image
of this embedding should lift.
It follows from Lemma~\ref{reduced} that
$\hat{\X}_t(\sigma(\O_{D}))$ is reduced.
Therefore it suffices to show that the generic
fibers of the irreducible components of
$\hat{\X}_t(\sigma(\O_{D}))$ all have
characteristic 0.  But this follows from the
fact that the reduction (mod $p$) of
a universal deformation of $G_0$ has
endomorphism ring $\Zp$ (cf.\
\cite[Th.\,1.1]{kk1}).

     Suppose $p\mid N^-$.  It follows from assumption
(\ref{assume1}) and Remark~\ref{split} that
$\O_{D}\otimes\Zp$ is the
ring of integers in a quadratic extension of $\Qp$.  Let
$(\hat{A},\hat{\imath})$ be a deformation of
$(\hat{A}_t,\hat{\imath}_t)$ such that
$\sigma(\O_{D})\otimes\Zp\subset
\End(\hat{A},\hat{\imath})$.  Since
$\End(\hat{A},\hat{\imath})$ is contained in
$\End(\hat{A}_t,\hat{\imath}_t)$, which is a local Eichler order
of type $(1,p^2)$, we must have $\O_D\otimes\Zp\cong\Zps$.
Therefore $\hat{A}$ has multiplication by
\begin{align}
(\ODel\otimes\Z_p)\otimes(\O_{D}\otimes\Zp)&\cong
\hat{\O}_{1,p}\otimes\Zps \\[.1cm]
&\cong\left\{\begin{bmatrix}a&b\\pc&d\end{bmatrix}:
a,b,c,d\in\Zps\right\}.
\end{align}
It follows that $\hat{A}\cong G\times G$,
where $G$ is a deformation of $G_0$.  Since
$\sigma(\O_D)\otimes\Zp$ commutes with the image of
$\hat{\imath}$, the map $\sigma$ induces an embedding
$\tau:\O_{D}\otimes\Zp\ra\End(G_0)$ such that
$\tau(\O_{D}\otimes\Zp)\subset\End(G)$.  Conversely,
any deformation $G$ of $G_0$ such that
$\tau(\O_{D}\otimes\Zp)\subset\End(G)$
gives a deformation $(\hat{A},\hat{\imath})$ of
$(\hat{A}_t,\hat{\imath}_t)$
such that $\sigma(\O_{D})\otimes\Zp
\subset\End(\hat{A},\hat{\imath})$.
The maximal deformation $G$ of $G_0$ such that
$\tau(\O_{D}\otimes\Zp)\subset\End(G)$
is the canonical lifting of $G_0$ associated to $\tau$.
By \cite[Prop.\,2.1]{qc} the canonical
lifting is defined over an integral domain
of characteristic 0.  As above this
implies that $\hat{\X}_t(\sigma(\O_{D}))$
is contained in the support of $\hat{\P}_{D,\pm b}$.
\\[\smallskipamount]
(b) Let $\hat{y}$ be an irreducible
component of $\hat{\P}_{D,\pm b}$, and let
$(A_y,i_y,Z_y)$ be the triple corresponding
to $y$.  Then there are two embeddings
$\rho_y,\rhobar_y:\O_D\ra\End(A_y,i_y,Z_y)$, where
$\rhobar_y$ is $\rho_y$ composed with complex conjugation on
$\O_D$.  Let
\begin{equation}
j:\End(A_y,i_y,Z_y)\lra\End(A_t,i_t,Z_t)
\end{equation}
be the natural embedding and define $\sigma=j\circ\rho_y$.
Then $\hat{y}$ is contained in $\X(\sigma(\O_D))$.  We need
to show that $\sigma:\O_D\ra\End(A_t,i_t,Z_t)$ is a
$\pm b$-embedding.

     Suppose $l\mid N$ with $l\not=p$.  Let $T_l(A_y)$,
$T_l(A_t)$ be the $l$-adic Tate modules of $A_y$, $A_t$, and
let $U_l^y\supset T_l(A_y)$, $U_l^t\supset T_l(A_t)$ be the
lattices which correspond to the $l$-primary parts of
$Z_y$, $Z_t$.  Then there is a natural
$(\ODelN\otimes\Zl)$-linear isomorphism $\nu:U_l^y\ra U_l^t$.
For $\gamma\in\End(A_y,i_y,Z_y)$ let $\tilde{\gamma}$ denote
the endomorphism of $U_l^y$ induced by $\gamma$, and let
$\widetilde{j(\gamma)}$ denote the endomorphism of $U_l^t$
induced by $j(\gamma)$.  Then we have
$\nu\circ\tilde{\gamma}=\widetilde{j(\gamma)}\circ\nu$.  It
now follows from the definitions of $\psi_l$ and $\omega_l$
that $\psi_l\circ j(\gamma)=\omega_l(\gamma)$.  Hence
\begin{equation} \label{psilsig}
\psi_l\circ\sigma=\psi_l\circ j\circ\rho_y=
\omega_l\circ\rho_y.
\end{equation}

     Suppose $p\mid N^-$.
It follows from Remark~\ref{split} that $p$ does not
split in $K$, and it follows from the fact that
$i_t$ has type 1 or 2 that $p$ is not ramified in
$K$.  Thus $p$ is inert in $K$.  Hence by
Proposition~\ref{product},
$T_p(A_y)/pT_p(A_y)$ is a vector space of dimension 2 over
$\R/p\R\cong\Fps$.  Furthermore, the
representation of $\End(A_y)$ on $T_p(A_y)/pT_p(A_y)$ is
equivalent to the representation of $\End(A_y)$ on
$\Lie(A_t)$ induced by $j$.  Therefore by the
constructions of $\psi_p$ and $\omega_p$ we get $\psi_p\circ
j\equiv\omega_p\pmod{p}$.  It follows that
\begin{equation} \label{psipsig}
\psi_p\circ\sigma\equiv\psi_p\circ j\circ\rho_y\equiv
\omega_p\circ\rho_y\pmod{p}.
\end{equation}
Since $y$ lies in the support of $P_{D,\pm b}$, it
follows from (\ref{psilsig}) and (\ref{psipsig}) that
$\sigma$ is a $\pm b$-embedding.
\\[\smallskipamount]
(c) Suppose $\hat{\X}_t(\sigma(\O_{D}))$ and
$\hat{\X}_t(\sigma'(\O_{D}))$ have a component
in common.  Then by (a) this component
is contained in the support of $\hat{\P}_{D,\pm b}$.
Let $y$ be the generic point of this component.
Then $y$ has characteristic 0, so by
Proposition~\ref{product}, $\End(A_y,i_y,Z_y)$
is an order in $\Q(\sqrt{D})$.  Since $\sigma(\O_{D})$ and
$\sigma'(\O_{D})$ both lie in $\End(A_y,i_y,Z_y)$, we
must have $\sigma(\O_{D})=\sigma'(\O_{D})$.~\qed

\begin{remark}
Suppose $p$ is ramified in $\Delta$
and $p$ divides the conductor $c$ of
$\O_{D}$.  Then $\hat{\X}_t(\sigma(\O_{D}))$
contains the subscheme of $\hat{\X}_t$
defined by the ideal $(p)$.  Therefore (a)
and (c) of Lemma~\ref{lift} are false if
$p\mid c$.  (However, (b) holds even if $p\mid c$.)
\end{remark}

\begin{prop} \label{combine}
We have $\hat{\P}_{D,\pm b}=
\frac12\cdot\sum\,\X_t(\sigma(\O_{D}))$,
where the sum is taken over all $\pm b$-embeddings
$\sigma:\O_{D}\ra\End(A_t,i_t,Z_t)$.
\end{prop}

\proof This follows from Lemma~\ref{lift}.
The factor $\frac12$ arises because $\sigma$ has the
same image as $\sigmabar$, where $\sigmabar$ is $\sigma$
composed with complex conjugation. \Qed

\section{Completion of the proofs}
\label{number}
\setcounter{equation}{0}

     In this section we use the results proved in
\S\ref{points}--\S\ref{universal} to compute
$\langle\P_{D_1,\pm b_1}\cdot\P_{D_2,\pm b_2}\rangle_{\X}$
and $\langle\P_{D_1}\cdot\P_{D_2}\rangle_{\X}$.  We start by
deriving a preliminary formula for $\langle\P_{D_1,\pm b_1}
\cdot\P_{D_2,\pm b_2}\rangle_{\X}$.

     It follows from (\ref{pdecomp}) that there are integers
$(\P_{D_1,\pm b_1}\cdot\P_{D_2,\pm b_2})_p$ such that
\begin{equation} \label{expressed}
\langle\P_{D_1,\pm b_1}\cdot\P_{D_2,\pm b_2}\rangle_{\X}=
\sum_{p<\infty}\;\;(\P_{D_1,\pm b_1}\cdot\P_{D_2,\pm b_2})_p
\cdot\log p.
\end{equation}
The quantity $(\P_{D_1,\pm b_1}\cdot\P_{D_2,\pm b_2})_p$
can be interpreted as an intersection multiplicity on
$\X\otimes W_p$.  Let $\I_p$ denote the set of points
$t\in\X(\Fpb)$ such that $(A_t,i_t,Z_t)\in\T_p$.  Then
the support of the intersection of
$\P_{D_1,\pm b_1}\otimes W_p$ with $\P_{D_2,\pm b_2}\otimes W_p$
on $\X\otimes W_p$ is contained in $\I_p$.  Therefore we have
\begin{equation} \label{prime}
(\P_{D_1,\pm b_1}\cdot\P_{D_2,\pm b_2})_p=
\sum_{t\in \I_p}\;\;(\P_{D_1,\pm b_1}\cdot\P_{D_2,\pm b_2})_{t}.
\end{equation}
Given embeddings $\sigma_1,\sigma_2$
of $\O_{D_1},\O_{D_2}$ into $\End(A_t,i_t,Z_t)$, set
$\epsilon_j=\sigma_j(\sqrt{D_j})$ and
$\gamma_j=(D_j+\epsilon_j)/2$ for $j=1,2$.  Then
$\gamma_j\in\End(A_t,i_t,Z_t)$ and
$\sigma_j(\O_{D_j})=\Zp[\gamma_j]$.  Hence by
Lemma~\ref{facts}(c) and Proposition~\ref{combine} we have
\begin{equation} \label{point}
(\P_{D_1,\pm b_1}\cdot\P_{D_2,\pm b_2})_{t}=
\frac{1}{4}\cdot\sum_{\sigma_1,\sigma_2}
\alpha_p(\gamma_1,\gamma_2),
\end{equation}
where the sum is taken over all pairs $(\sigma_1,\sigma_2)$
such that $\sigma_j$ is a $\pm b_j$-embedding.

     It follows from Proposition~\ref{orders} that
$\End(A_t,i_t,Z_t)$
is isomorphic to an Eichler order of type $(N^+,pN^-)$
in the quaternion algebra $\Delta(p)$ over $\Q$.  Therefore
$\End(A_t,i_t,Z_t)$ has reduced discriminant $N^+\cdot pN^-=pN$.
Let $n=\frac12\cdot\Tr(\epsilon_1\epsilon_2)=
2\Tr(\gamma_1\gamma_2)-D_1D_2$, where $\Tr$ is the reduced
trace from $\Delta(p)$ to $\Q$. Since $\sigma_j$ is a
$\pm b_j$-embedding we have $n\equiv
\pm h\pmod{2N}$, where $h$ is determined by (\ref{h1}) and
(\ref{h2}).  Since $n\equiv D_1D_2\pmod2$, the ring
$S_n$ from \S\ref{eichler} is defined.  In fact
$S_n\cong\Z[\gamma_1,\gamma_2]$ is isomorphic to a suborder
of $\End(A_t,i_t,Z_t)$, and hence $B_n\cong\Delta(p)$ and
$pN$ divides the reduced discriminant
$\delta_n=(n^2-D_1D_2)/4$ of $S_n$.
Since $\Delta(p)$ is ramified at $\infty$
we have $\delta_n<0$, and hence $n^2<D_1D_2$. 
It follows that $p\le|\delta_n|/N\le D_1D_2/4N$, so by
assumption C3 in (\ref{cond}) we have $p\nmid m$.

     The ternary quadratic form $\Nr(x+y\gamma_1+z\gamma_2)$
over $\Z$ is equal to the form $Q_n(x,y,z)$ defined in
(\ref{Qnexp}).  In \S\ref{universal} we saw that the
intersection number $\alpha_p(\gamma_1,\gamma_2)$ depends
only on the $\Zp$-isometry class of $Q_n$ and not on the
choice of $\gamma_1,\gamma_2$. Therefore we may write
$\alpha_p(Q_n)=\alpha_p(\gamma_1,\gamma_2)$.
Using (\ref{prime}) and (\ref{point}) we get
\begin{align}
(\P_{D_1,\pm b_1}\cdot\P_{D_2,\pm b_2})_p&=
\sum_{t\in \I_p}\;\;\left(\frac{1}{4}\cdot
\sum_{\sigma_1,\sigma_2}\alpha_p(\gamma_1,\gamma_2)
\right) \\[.1cm]
&=\frac{1}{4}\cdot\underset{n\equiv\pm h\;(2N)}{\sum_{n^2<D_1D_2}}\;\left(\sum_{t\in \I_p}\;
r_{t}(n,\pm b_1,\pm b_2)\right)\cdot
\alpha_p(Q_n), \label{formula}
\end{align}
where $r_{t}(n,\pm b_1,\pm b_2)$ denotes the number of pairs
$(\sigma_1,\sigma_2)$ of $\pm b_j$-embeddings of
$\O_{D_1},\O_{D_2}$ into $\End(A_t,i_t,Z_t)$ such that
$\frac{1}{2}\cdot\Tr(\epsilon_1\epsilon_2)=n$.

     We wish to interpret the inner sum of (\ref{formula}) as
counting embeddings of the ring $S_n$ into Eichler orders.
Let $\E\subset\Delta(p)$ be an Eichler order of type
$(N^+,pN^-)$, and let
\begin{equation}
N(\E)=\{\beta\in\Delta(p)^{\times}:\beta \E\beta^{-1}=\E\}.
\end{equation}
Then $N(\E)/\Q^{\times}\E^{\times}$ acts freely on the set of
orientations on $\E$; the orbits of this action are the
isomorphism classes of orientations for $\E$.
Let $v_{\E}=|N(\E)/\Q^{\times}\E^{\times}|$ and recall that
$r$ is the number of distinct primes dividing $N$.
Thus if $p\nmid N^-$ then $\E$ has $2^{r+1}$ orientations,
and hence $2^{r+1}/v_{\E}$ isomorphism classes of
orientations, while if
$p\mid N^-$ then $\E$ has $2^r$ orientations, and hence
$2^r/v_{\E}$ isomorphism classes of orientations.
By Proposition~\ref{orders}, $\End(A_t,i_t,Z_t)$ is an
Eichler order of type $(N^+,pN^-)$ for each $t\in\I_p$.
Let $\I_{\E}$ denote the set of $t\in\I_p$ such that
$\End(A_t,i_t,Z_t)\cong \E$.  Let $\c_p$ be a set of
representatives for the isomorphism classes of Eichler orders
of type $(N^+,pN^-)$.  Then
\begin{equation} \label{rtnut}
\sum_{t\in\I_p}r_t(n,\pm b_1,\pm b_2)=
\sum_{\E\in\c_p}\sum_{t\in\I_{\E}}r_t(n,\pm b_1,\pm b_2).
\end{equation}

     Given an orientation $\{\mu_l\}_{l\mid pN}$ on
the Eichler order $\E$, say that the homomorphism
$\tau:S_n\ra \E$ is a $(b_1,b_2)$-embedding if
$\mu_l\circ\tau(e_j)\equiv(b_j)_l\pmod{l^{a_l}}$ for $j=1,2$
and all $l\mid N$, where $a_l=v_l(2N)$.  Say that $\tau$ is a
$(\pm b_1,\pm b_2)$-embedding if $\tau$ is a $(b_1,b_2)$-,
$(-b_1,b_2)$-, $(b_1,-b_2)$-, or $(-b_1,-b_2)$-embedding.  Let
\begin{equation}
\T_{\E}=\{(A,i,Z)\in\T_p:\End(A,i,Z)\cong\E\}
\end{equation}
and let
$\tilde{\I}_{\E}$ be a subset of $\I_{\E}$ such that each
$(A,i,Z)\in\T_{\E}$ is isomorphic to $(A_t,i_t,Z_t)$ for
exactly one $t\in\tilde{\I}_{\E}$.  

     Let $r_{\E}(n)$ denote
the number of embeddings of $S_n$ into $\E$.  Suppose
$p\nmid N^-$.  Then each homomorphism $\tau:S_n\ra\E$ is a
$(\pm b_1,\pm b_2)$-embedding for four different orientations
on $\E$.  Hence the set
\begin{equation}
\Upsilon_n=\{(\tau:S_n\ra\E,\{\mu_l\}_{l\mid pN}):
\text{$\tau$ is a $(\pm b_1,\pm b_2)$-embedding w.\,r.\,t.\ 
$\{\mu_l\}_{l\mid pN}$}\}
\end{equation}
has cardinality $4r_{\E}(n)$.  It follows from
Proposition~\ref{crossing}(a) that the isomorphism classes
of orientations on $\E$ correspond to elements
$t\in\tilde{\I}_{\E}$.  This allows us to count the elements
of $\Upsilon_n$ in a different manner, and gives the
formula
\begin{equation} \label{vRrt}
4r_{\E}(n)=
v_{\E}\cdot\sum_{t\in\tilde{\I}_{\E}}r_t(n,\pm b_1,\pm b_2).
\end{equation}
Suppose $p\mid N$.  Then by
Proposition~\ref{crossing}(b) each isomorphism class of
orientations on $\E$ is represented by two different
$t\in\tilde{\I}_{\E}$, and each homomorphism
$\tau:S_n\ra \E$ is a $(\pm b_1,\pm b_2)$-embedding for
two orientations on $\E$.  Hence (\ref{vRrt}) is valid in
this case as well. 

     Let $(A,i,Z)\in\T_{\E}$.  Then the group
$\Aut(A,i,Z)\cong\E^{\times}$ acts freely on the set of
level-$m$ structures $\beta$ on $(A,i)$.  By
Remark~\ref{eta} the pair $(A,i)$ admits $2\eta(m)$
different level-$m$ structures.  Hence there are $2\eta(m)/2u_{\E}$
isomorphism classes of 4-tuples $(A,i,Z,\beta)$,
where $u_{\E}=\frac12\cdot\#\E^{\times}$.  It follows that
\begin{equation}
\sum_{t\in\I_{\E}}r_t(n,\pm b_1,\pm b_2)=
\frac{\eta(m)}{u_{\E}}\cdot
\sum_{t\in\tilde{\I}_{\E}}r_t(n,\pm b_1,\pm b_2).
\end{equation}
Combining this formula with (\ref{vRrt}) we get
\begin{equation} \label{etare}
\sum_{t\in\I_{\E}}r_t(n,\pm b_1,\pm b_2)=
\frac{4\eta(m) r_{\E}(n)}{u_{\E}v_{\E}}.
\end{equation}

     It follows from (\ref{rtnut}) and (\ref{etare}) that
\begin{align} \label{uRvR}
\sum_{t\in\I_p}r_t(n,\pm b_1,\pm b_2)
&=4\eta(m)\cdot\sum_{\E\in\c_p}\frac{r_{\E}(n)}{u_{\E}v_{\E}}.
\end{align}
By combining this formula with (\ref{formula})
and (\ref{expressed}) we get the following formula
for the arithmetic intersection number:
\begin{equation} \label{prelim}
\langle\P_{D_1,\pm b_1}\cdot\P_{D_2,\pm b_2}\rangle_{\X}=
\eta(m)\cdot\sum_{p<\infty}
\left(\underset{n\equiv\pm h\;(2N)}{\sum_{n^2<D_1D_2}}\left(
\sum_{\E\in \c_p}\;\frac{r_{\E}(n)}{u_{\E}v_{\E}}
\right)\cdot\alpha_p(Q_n)\right)\cdot\log p.
\end{equation}
To prove Theorem~\ref{repnum} we will evaluate the inner
sum of (\ref{prelim}) in terms of representation numbers of
quadratic forms.
To prove Theorem~\ref{explicit} we will evaluate the inner
sum of (\ref{prelim}) by counting Eichler orders. \medskip

\noindent {\em Proof of Theorem~\ref{repnum}:}
Recall that $r$ denotes the number of distinct prime divisors
of $N$.  For each $n$ such that $n^2\equiv D_1D_2\pmod{2N}$
there are $2^{r-1}$ pairs $(\pm b_1,\pm b_2)$ such that
$n\equiv\pm h(\pm b_1,\pm b_2)\pmod{2N}$.  If
$r_{\E}(n)\not=0$ then $pN$
divides the reduced discriminant $\delta_n=(n^2-D_1D_2)/4$ of
$S_n$, so we have $n^2\equiv D_1D_2\pmod{4pN}$.
Therefore by summing (\ref{prelim}) over all $2^{2r-2}$ pairs
$(\pm b_1,\pm b_2)$ we get
\begin{equation} \label{prelim2}
\langle\P_{D_1}\cdot\P_{D_2}\rangle_{\X}=
2^{r-1}\cdot\eta(m)\cdot\sum_{p<\infty}
\left(\underset{n^2\equiv D_1D_2\;(4pN)}{\sum_{n^2<D_1D_2}}\left(
\sum_{\E\in\c_p}\;\frac{r_{\E}(n)}{u_{\E}v_{\E}}
\right)\cdot\alpha_p(Q_n)\right)\cdot\log p.
\end{equation}

     Let $\E\in\c_p$ and let $n$ be an integer such that
$n^2<D_1D_2$ and $n^2\equiv D_1D_2\pmod{4pN}$.  Recall
that $S_n=\Z[g_1,g_2]$ and define
$L_n=\Z+\Z g_1+\Z g_2\subset S_n$.  Then $L_n$ is a quadratic
space with the reduced norm form $\Nr$.
Let $(\tau,\lambda)$ be a pair consisting of a ring
homomorphism $\tau:S_n\ra \E$ and a unit
$\lambda\in \E^{\times}$.  Associated to this pair is an
isometry $\lambda\cdot\tau|_{L_n}$ of $L_n$ into $\E$.  We
claim that every isometry from $L_n$ into $\E$ arises this
way.  Let $\nu:L_n\ra \E$ be an isometry, and set
$\lambda=\nu(1)$.  Then $1=\Nr(1)=\Nr(\lambda)$, so
$\lambda\in \E^{\times}$.  Define a $\Z$-linear map
Let $\tau:S_n\ra\E$ be the unique $\Z$-linear map such that
$\tau(1)=1$, $\tau(g_j)=\lambda^{-1}\nu(g_j)$, and
$\tau(g_1g_2)=\tau(g_1)\tau(g_2)$.
Set $\epsilon_j=\tau(e_j)$ for $j=1,2$.  Since $\nu$ is an
isometry we get $\epsilon_j^2=-\Nr(\epsilon_j)=-\Nr(e_j)=D_j$
and
\begin{align}
(\epsilon_1+\epsilon_2)^2&=-\Nr(\epsilon_1+\epsilon_2) \\
&=-\Nr(e_1+e_2) \\
&=(e_1+e_2)^2,
\end{align}
which implies
$\epsilon_1\epsilon_2+\epsilon_2\epsilon_1=2n$.  It follows
that
\begin{equation}
\tau\otimes\Q:S_n\otimes\Q\lra\E\otimes\Q
\end{equation}
is an isomorphism of rings, and hence that $\tau$ is a ring
homomorphism.

     Since $\Z^3$ with the quadratic form $Q_n$ is isometric
to $L_n$, it follows from the preceding paragraph that 
$\R_{\E}(Q_n)=\#\E^{\times}\cdot r_{\E}(n)=2u_{\E}r_{\E}(n)$.
By Lemma~\ref{self} the number of
proper self-isometries of $\E$
is $w_{\E}=2u_{\E}^2v_{\E}$.  Therefore
$r_{\E}(n)/u_{\E}v_{\E}=\R_{\E}(Q_n)/w_{\E}$.  Hence
by Proposition~\ref{S} we have
\begin{equation}
\sum_{\E\in\c_p}\frac{r_{\E}(n)}{u_{\E}v_{\E}}=
\sum_{L\in\G_p}\frac{\R_L(Q_n)}{w_L}.
\end{equation}
Substituting this formula into (\ref{prelim2}) gives
Theorem~\ref{repnum}. \qed

\noindent {\em Proof of Theorem~\ref{explicit}:}
Choose $n$ such that $n^2<D_1D_2$ and $n\equiv\pm h\pmod{2N}$.
Since $\gcd(D_1,D_2)=1$, it
follows from \S\ref{eichler} that $S_n$ is an Eichler order
of type $(\delta_n^+,\delta_n^-)$, and it follows from
Remark~\ref{Ndel} that $N^+\mid \delta_n^+$ and
$N^-\mid \delta_n^-$.  Set $M_n^+=\delta_n^+/N^+$,
$M_n^-=\delta_n^-/N^-$,
$m_p^+=v_p(M_n^+)$, and $m_p^-=v_p(M_n^-)$, and define
$L_{(p)}(s)=L_{p^{m_p^+},\:p^{m_p^-}}(s)$.  Then we have
\begin{align}
L_{M_n^+,M_n^-}(s)&=\prod_{p<\infty}L_{(p)}(s) \\
L_{M_n^+,M_n^-}'(0)&=\sum_{p<\infty}
\left(\prod_{l\not=p}L_{(l)}(0)\right)L_{(p)}'(0).
\label{der2}
\end{align}

     Let $p$ be a prime, and assume for now that
$S_n\otimes\Q\cong\Delta(p)$.  If $p\mid N^+$ then
$p\nmid\delta_n^-$ by Remark~\ref{Ndel}, so this assumption
implies $p\nmid N^+$.  Let $\E_0$ be an Eichler order
of type $(N^+,pN^-)$.  We may compute $r_{\E_0}(n)$ as the
number of embeddings of $\E_0$ into $S_n\otimes\Q$ whose
image contains $S_n$.  Hence $r_{\E_0}(n)$ is the sum over all
$\E$ such that
${S_n\subset\E\subset S_n\otimes\Q}$ of the number of
isomorphisms of $\E_0$ onto $\E$.  If $\E_0\cong\E$ then
there are $|N(\E_0)/\Q^{\times}|=u_{\E_0}v_{\E_0}$ such
isomorphisms.  We conclude that the number of Eichler orders
$\E$ such that $S_n\subset\E\subset S_n\otimes\Q$ and
$\E\cong\E_0$ is $r_{\E_0}(n)/u_{\E_0}v_{\E_0}$.  It follows
that the inner sum $\sum_{\E\in\c_p}r_{\E}(n)/u_{\E}v_{\E}$
of (\ref{prelim}) is
equal to the number of Eichler orders $\E$ of type
$(N^+,pN^-)$ such that $S_n\subset \E\subset S_n\otimes\Q$.

     For each prime $l$ set $e_l^+=v_l(N^+)$, $e_l^-=v_l(pN^-)$,
$d_l^+=v_l(\delta_n^+)$, and
$d_l^-=v_l(\delta_n^-)$. The number of Eichler orders in
$ S_n\otimes\Q$ of type $(N^+,pN^-)$ which contain the Eichler
order $S_n$ can be computed as the product over $l$ of the
number $c_l$ of local Eichler orders
$\E_l$ of type $(l^{e_l^+},l^{e_l^-})$ in $S_n\otimes\Q_l$
such that $\E_l\supset S_n\otimes\Zl$.  
If $S_n\otimes\Zl$ has type $(l^{d_l^+},1)$ then
$c_l=d_l^+-e_l^++1=L_{(l)}(0)$.  If $S_n\otimes\Zl$
has type $(1,l^{d_l^-})$ then $c_l=1$ for
$d_l^-\equiv e_l^-\pmod2$ and $c_l=0$ for
$d_l^-\not\equiv e_l^-\pmod2$, so once again we have
$c_l=L_{(l)}(0)$.  Since we are assuming
$S_n\otimes\Q\cong\Delta(p)$ we get $c_p=1$.
Hence the total number of Eichler orders $\E$ of type
$(N^+,pN^-)$ in $\Delta(p)$ which contain $S_n$ is
\begin{equation} \label{LAN}
\sum_{\E\in\c_p}\frac{r_{\E}(n)}{u_{\E}v_{\E}}=
\prod_{l\not=p}L_{(l)}(0).
\end{equation}
Now suppose that $S_n\otimes\Q$ is not isomorphic to
$\Delta(p)$.  In this case the sum on the left side of
(\ref{LAN}) is 0, and there
is at least one prime $l\not=p$ such that $v_l(M_n^-)$ is
odd.  It follows that $L_{(l)}(0)=0$, so the right side of
(\ref{LAN}) is also 0.
Hence (\ref{LAN}) is valid for all primes $p$. 

     Suppose $p\nmid N$.  Then the assumption $\gcd(D_1,D_2)=1$
implies that, in the notation of (\ref{unram}), we have
$a_2=0$ and $a_3=v_p(\delta_n)=v_p(\delta_n^-)$.  Therefore
$\alpha_p(Q_n)=(v_p(\delta_n^-)+1)/2$.
On the other hand, if $p\mid N^-$ then by (\ref{ram}) we have
$\alpha_p(Q_n)=v_p(\delta_n^-)/2$.
In both cases we get $\alpha_p(Q_n)=(v_p(M_n^-)+1)/2$.
If $S_n\otimes\Q\cong\Delta(p)$ then $v_p(M_n^-)$ is odd.
Therefore by (\ref{L1pe}) we get
\begin{align}
L_{(p)}'(0)&=\frac{v_p(M_n^-)+1}{2}\cdot\log p \\
&=\alpha_p(Q_n)\cdot\log p. \label{alog}
\end{align}
It follows from (\ref{der2}), (\ref{LAN}), and (\ref{alog}) that
\begin{equation}
L_{M_n^+,M_n^-}'(0)=\sum_{p<\infty}
\left(\sum_{\E\in\c_p}\frac{r_{\E}(n)}{u_{\E}v_{\E}}\right)
\cdot\alpha_p(Q_n)\cdot\log p.
\end{equation}
Combining this formula with (\ref{prelim}) we get
\begin{equation} \label{final}
\langle\P_{D_1,\pm b_1}\cdot\P_{D_2,\pm b_2}\rangle_{\X}=
\eta(m)\cdot
\underset{n\equiv\pm h\;(2N)}{\sum_{n^2<D_1D_2}}
\;L_{M_n^+,M_n^-}'(0),
\end{equation}
which is Theorem~\ref{explicit}. \Qed

\end{document}